\documentclass[12pt]{article}

\usepackage[dvips]{graphicx}
\usepackage[latin1]{inputenc}
\usepackage{psfrag}
\usepackage{amssymb}

\newcounter{sec}

\newcommand{\eq}[1]{$$ #1 $$}
\newcommand{\eqa}[1]{\begin{eqnarray} #1 \end{eqnarray}}
\newcommand{\n}{\noindent}
\newcommand{\ite}[3]{\n
\begin{minipage}[c]{\textwidth}\vspace{0.4cm}{\bf #1 \thesection.\thesec\ #2\,} \emph{#3} \addtocounter{sec}{1}\vspace{0.3cm}\end{minipage}}
\newcommand{\proof}{\n\emph{Proof\,. }}
\newcommand{\qed}{\hspace{0.5cm} $\square$\\}
\newcommand{\h}{\mathbf{H}}
\newcommand{\s}{\mathbf{S}}
\newcommand{\e}{\mathbf{E}}
\newcommand{\R}{\mathbf{R}}
\newcommand{\z}{\mathbb{Z}}
\newcommand{\esc}[1]{\langle #1 \rangle}

\begin{document}
\thispagestyle{empty}
\begin{center}
\hspace*{10cm} {\bf  } \\
\hspace*{8cm} {\bf   }\\
\begin{Large}
{\bf On the Characterization of Polyhedra in\\ Hyperbolic 3-Space}
 \\[1cm]
\end{Large}
{ \sc Javier Virto}  \\[0.2cm]
\begin{center}
{\em 
INFN, Sezione di Roma, I-00185 Rome, Italy
}
\end{center}



\vspace{0.5cm}
%
\end{center}

\tableofcontents

\vspace{0.1cm}

\section{Introduction}

\setcounter{sec}{1}

The issue of the characterization of possible polytopes in 3 dimensional space is an old subject. There are some very interesting theorems that were proven long ago. Probably one of the most famous is Cauchy's Rigidity Theorem, which states that an euclidean polytope with rigid faces cannot be deformed. The condition of rigid faces cannot be softened; indeed, face angles do not determine an euclidean polyhedron, and neither do edge lengths, as can be seen easily by deforming cubes.

When one moves to other non-euclidean geometries, things get often more interesting. There is a property of homogeneously curved geometries: the existence of a length \emph{scale} that breaks scale invariance. For example, all regular right-angled octagons in the hyperbolic plane have the same size: if you want a bigger one, your angles will decrease. A more familiar example is a right-angled spherical triangle. Therefore one should not be surprised by the fact that indeed hyperbolic and spherical polyhedra \emph{are} determined by face angles, although they are \emph{not} determined by edge lengths. All these issues are worth studying.

Concerning hyperbolic polyhedra, much is already known. Convex hyperbolic polyhedra with all vertices trivalent are determined by their dihedral angles, as well as ideal polyhedra (those with the vertices at infinity). Andreev gave a complete characterization of compact convex hyperbolic polyhedra with dihedral angles not larger than $\pi /2$, and Rivin and coworkers have completed such characterizations for general compact, ideal and finite volume, convex hyperbolic polyhedra.

The purpose of these pages is to review several results related to such characterizations of polyhedra in hyperbolic 3-space. In particular we present Rivin's theorem that gives a characterization of compact convex hyperbolic polyhedra, and Hodgson's proof of the Adreev's theorem. We also review the analogous characterization of ideal polyhedra, and give a family of counter-examples that proves that hyperbolic polyhedra are not determined by edge lengths.

\section{Definitions, models, constructions and dualities}

\setcounter{sec}{1} \label{defs}

We begin with some basic definitions and constructions. These are standard and can be found, for example, in the books by Thurston \cite{Th} and
Ratcliffe \cite{Rat}, and in the paper by Rivin and Hodgson \cite{Hodgson}.

Euclidean, Hyperbolic and Spherical n-spaces are the $n$-dimensional simply connected Riemannian manifolds with constant curvature $0$, $-1$ and $1$
respectively. They are denoted by $\e^n$, $\h^n$ and $\s^n$. Lorentzian, de Sitter, and anti de Sitter $n$-spaces are the symmetric pseudo-Riemannian
manifolds of signature $(n-1,1)$ (Lorentzian manifolds), with constant curvature $0$, $1$ and $-1$ respectively. They are denoted by $\e^{n-1}_1$,
$\s^{n-1}_1$ and $\h^{n-1}_1$. The basic models for $\e^n$ and $\e^{n-1}_1$ are $\R^n$ with the Euclidean and Minkowski inner products. The basic
model for $\s^n$ is the $n$-sphere embedded in $\e^{n+1}$ with the induced metric. The de Sitter $n$-space can be modeled on a hyperboloid of one sheet
in $\e^{n}_1$:
\eq{\s^{n-1}_1=\{x\in \e^{n}_1\,|\, \esc{x,x}=1\},}
and it can be easily seen that the induced metric is Lorentzian. In the same way the anti de Sitter $n$-space can be modeled as
\eq{\h^{n-1}_1=\{x\in \e^{n-1}_2\,|\, \esc{x,x}=-1\}.}
Here we are assuming $n\ge 2$ (we will work with $n=3$ specifically) and the convention used for Minkowski's metric is
$\esc{x,x}=-x_0^2+x_1^2+\cdots+x_{n-1}^2=-x_0^2+\|\vec{x}\|^2$. Also, a \emph{hyperplane} in any of these $n$-spaces is defined as a
$(n-1)$-dimensional geodesic plane.

We shall be more concerned about different models for $\h^n$. The \emph{hyperboloid model} for $\h^n$ is the upper sheet of the two-sheeted
hyperboloid in $\e^{n}_1$,
\eq{\h^{n}=\{x\in \e^{n}_1\,|\, x_0>0,\,\esc{x,x}=-1\}.}
It is straightforward to see that the induced metric is indeed Riemannian and of constant curvature $-1$. The geodesic $m$-planes in $\h^n$ are the
intersection with the hyperboloid of $(m+1)$-planes in $\e^n_1$ passing through the origin.

From the hyperboloid model arises a projective model, the \emph{Klein model} or just \emph{projective model} for $\h^n$. Each point $v$ on the
hyperboloid is mapped to a point in the unit disk $D^n=\{x\in \e^n_1\,|\,x_0=1,\,\|\vec{x}\|^2< 1\}$: the point where the line joining $v$ with the
origin of $\e^n_1$ intersects $D^n$. This model is projective, that is, hyperbolic $m$-planes are the intersection of Euclidean $(m+1)$-planes with
$D^n$ (which are Euclidean $m$-disks), but it is not conformal, since angles do not correspond to Euclidean angles.

A conformal model on the unit $n$-disk is provided by the \emph{Poincare model} or \emph{disk model} of $\h^n$. This is $D^n$ together with the
metric
\eq{ds^2=\frac{4}{(1-r^2)^2}(dx_1^2+\cdots+dx_n^2)}
where $r$ is the radial coordinate. Hyperbolic $m$-planes in this model are $m$-spheres intersecting orthogonally with $\partial D^n$, called the
\emph{sphere at infinity}, including Euclidean $m$-planes passing through the origin of the disk. This model is conformal since angles correspond to
Euclidean angles.

Finally, another useful model for $\h^n$ is the \emph{upper half-space model}. This is the space $\R^n_+=\{x\in \R^n\,|\,x_n\ge 0\}$ together with
the metric
\eq{ds^2=\frac{1}{x_n^2}(dx_1^2+\cdots+dx_n^2).}
Hyperbolic $m$-planes are given in this model by Euclidean $m$-planes and $m$-spheres orthogonal to the plane $x_n=0$.

\begin{figure}
\begin{center}
\psfrag{a}{(a)}\psfrag{b}{(b)}\psfrag{p}{\tiny$\pi/2$}
\includegraphics[height=7cm,width=5.5cm]{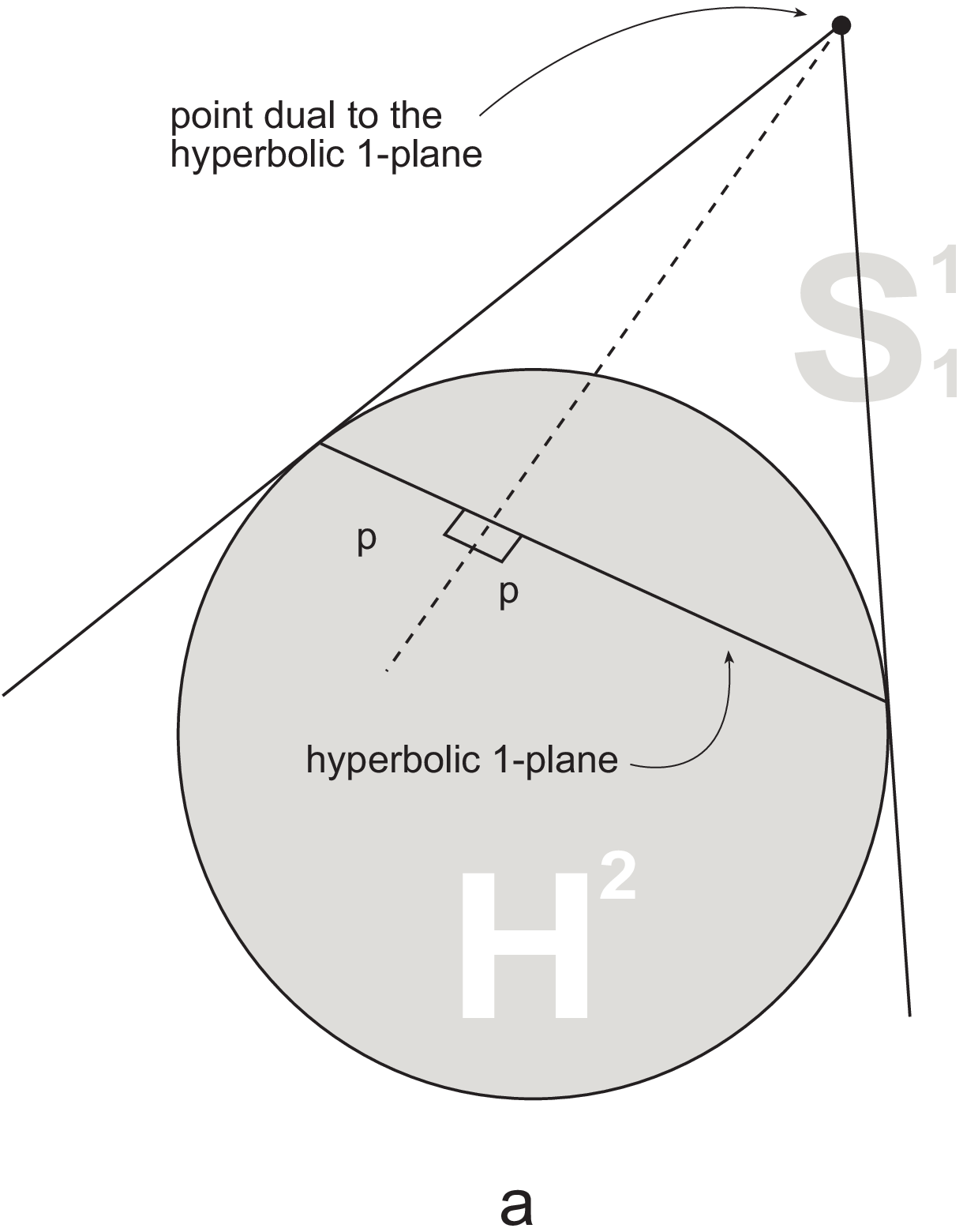}\hspace{2cm}\includegraphics[height=7cm,width=6.5cm]{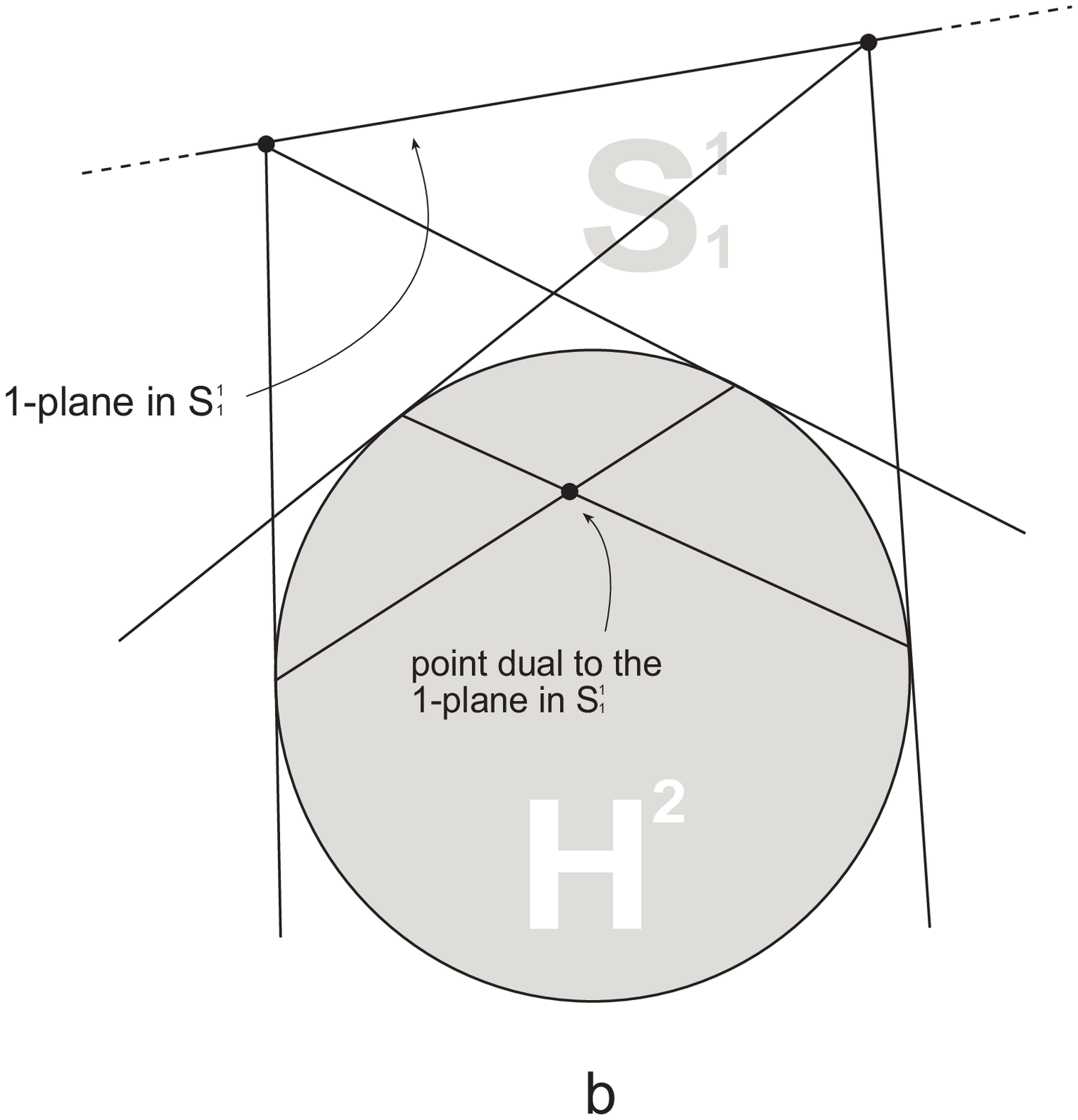}
\end{center}
\caption{(a) A $1$-plane in $\h^2$ can be represented by a point in $\s^1_1$ by means of the projective model. All the (Euclidean) lines passing
through that point and intersecting the hyperbolic $1$-plane are perpendicular to it. (b) Analogously, a $1$-plane in $\s^1_1$ can be represented by
a point in $\h^2$.} \label{f2.1}
\end{figure}

Oriented hyperplanes in $\e^n$ passing through the origin are specified by the unit vector orthogonal to the hyperplane, which is unique as specified
by the orientation. Therefore oriented hyperplanes in $\e^n$ are parameterized by the points on $\s^{n-1}$. In the same way, oriented hyperplanes in
$\s^n$ are parameterized by points on $\s^n$ (this is the usual \emph{polar map}). Oriented hyperplanes in $\h^n$ correspond to oriented
\emph{time-like} $n$-planes in $\e^{n}_1$ passing through the origin, so their normal vectors are space-like, and therefore these hyperplanes are
parameterized by points on $\s^{n-1}_1$. Similarly, oriented hyperplanes in $\s^{n-1}_1$ are parameterized by points on $\h^n$. This is the basics of
the \emph{polar duality}, to be defined later. A pictorial example of how planes in $\h^2$ and $\s^1_1$ are parameterized by points in $\s^1_1$ and
$\h^2$ respectively is shown in Fig.~\ref{f2.1} using the projective model.

The projective model is particularly suitable to visualize the polar duality between $\h^n$ and $\s^{n-1}_1$. In fact, in an analogous way to the
projective model of $\h^n$ one can model $\s^{n-1}_1$ projectively on $\hat{D}^n:=\{x\in \e^n_1\,|\,x_0=1,\,\|\vec{x}\|^2>1\}$ by sending a point $v$
on $\s^{n-1}_1$ (inside $\e^n_1$) to the intersection of the hyperplane $x_0=1$ with the line joining $v$ and the origin of $\e^n_1$. Then both
$\h^n$ and $\s^{n-1}_1$ are modeled together inside $\mathbf{RP}^n$. The hyperplane in $\h^n$ dual to a point $v$ in $\s^{n-1}_1$ is the Euclidean
hyperplane containing the sphere $\s^{n-2}$ of tangency of $\partial D^n$ with a hypercone with vertex at $v$. The hyperplane in $\s^{n-1}_1$ dual to
a point $v$ in $\h^n$ is the collection of points dual to hyperplanes in $\h^n$ that contain $v$. Duals of $m$-planes for any $m$ are found
analogously.

A \emph{convex polyhedron} is defined as the intersection of a finite number of half-spaces. If it is bounded, then it is just a traditional
polytope, but a convex polyhedron might be unbounded. A \emph{general polyhedron} is the union of a finite number of convex polyhedra. A
\emph{compact} polyhedron in $\h^3$ is a bounded hyperbolic polyhedron with all its vertices finite. An \emph{ideal} polyhedron in $\h^3$ is the
convex hull of a set of points located on the sphere at infinity, called \emph{ideal vertices}. An ideal polyhedron is not compact, but it has finite
volume (note, for example, that an ideal hyperbolic triangle in $\h^2$ has an area equal to $\pi$). A \emph{finite-volume} polyhedron in $\h^3$ is a
hyperbolic polyhedron with finite and ideal vertices. A \emph{hyperinfinite vertex} is a vertex beyond the sphere at infinity in the projective model
of hyperbolic space. Such a vertex is the dual of a hyperbolic hyperplane orthogonal to all the faces that (when extended) meet the vertex, so that a
a face of a hyperbolic polyhedron with all internal dihedral angles at its edges equal to $\pi/2$, can be represented by a hyperinfinite vertex.

\begin{figure}
\begin{center}
\psfrag{a}{(a)}\psfrag{b}{(b)}
\includegraphics[height=6cm,width=13cm]{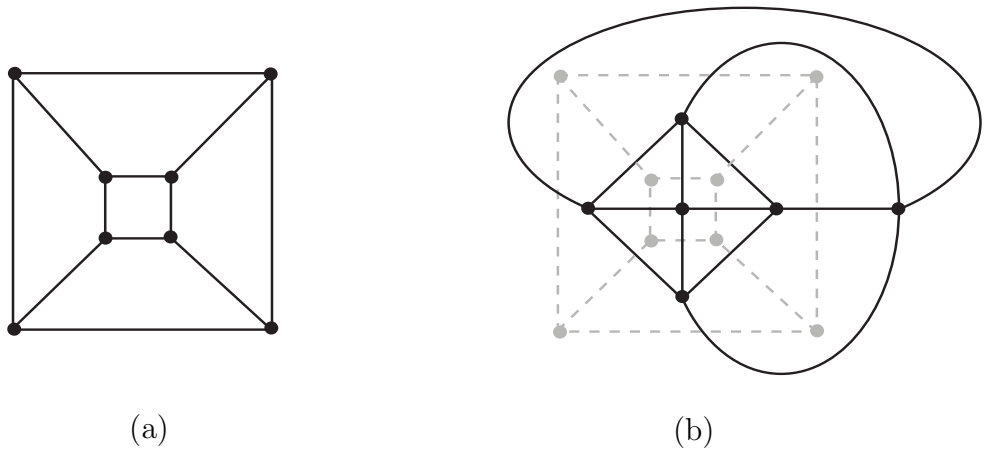}
\end{center}
\caption{(a) An abstract cube, represented by the 3-connected planar graph associated with the 1-skeleton of an Euclidean cube. (b) The Poincaré dual
of a cube, which is clearly an octahedron.} \label{f2.2}
\end{figure}

An \emph{abstract polyhedron} is a partially ordered set of elements, which in three dimensions fall into three categories that can be called
conveniently \emph{vertices}, \emph{edges} and \emph{faces}. Therefore the defining properties of an abstract polyhedron are the number of vertices,
edges and faces and its \emph{combinatorics}. More specifically, let $X_0$, $X_1$ and $X_2$ be the sets of oriented vertices, edges and faces
respectively, and define the $n$-chains $\Delta_n$ as the free abelian groups with basis $X_n$. Then the abstract polyhedron $P$ is defined by the
chain complex
\eq{0\stackrel{d}{\longrightarrow}\Delta_2\stackrel{d}{\longrightarrow}\Delta_1\stackrel{d}{\longrightarrow}\Delta_0\to 0}
where $d(f)=e_1+\cdots+e_n$ if $e_1,\dots,e_n$ bound the face $f$ cyclically, and $d(e)=v_2-v_1$ if the edge $e$ goes from $v_1$ to $v_2$ in this
direction. One can easily verify that indeed $d^2=0$. The \emph{Poincaré dual} of $P$ is defined by the dual cochain complex with chain groups
$\Delta_{n-k}^*=Hom(\Delta_{n-k},\z)$ and isomorphisms $\Delta_n\to\Delta_{n-k}^*$. Under this duality, faces are dual to vertices, and edges are
dual to edges. Fig.~\ref{f2.2} shows an abstract cube and its Poincaré dual.

A combinatorial characterization of convex polyhedra is given by a theorem due to Steinitz:

\ite{Theorem}{(Steinitz)}{A graph is the one-skeleton of a convex polyhedron in $\e^3$ if and only if it is a 3-connected planar graph.}

A $k$-connected planar graph is a planar graph such that after removing any $k-1$ edges is still connected. The Poincaré dual of an Euclidean
polyhedron has a 1-skeleton which is also a 3-connected planar graph, so Poincaré duals to Euclidean polyhedra are also Euclidean polyhedra. (For
example, in Fig.~\ref{f2.2} it is shown that the Poincaré dual of a cube is an octahedron.) Moreover, as can be easily inferred using the projective
model of $\h^3$, ideal hyperbolic polyhedra are in one to one correspondence with Euclidean convex polyhedra inscribed in the sphere, so ideal
polyhedra and their Poincaré duals have 1-skeletons which are 3-connected planar graphs. In addition, a graph specifies an abstract polyhedron and
its dual.

Let $P$ be a convex polyhedron in 3-space. Consider a vertex $v$ of the polyhedron, and place a small sphere centered at $v$. The intersection of $P$
with this sphere is a convex spherical polygon. The \emph{link} of the vertex $v$, $lk(v)$, is defined as this spherical polygon after rescaling the
sphere so that its radius is one. Now consider a convex spherical polygon $p$. The \emph{spherical polar} of $p$, denoted $p^*$, is the convex
spherical polygon whose vertices are the centers of the exterior hemispheres tangent to the faces of $p$.

\ite{Definition}{}{The {\bf generalized Gauss Map} $G$ is an application from the set of convex polyhedra in 3-space to the set of 2-dimensional
metric spaces, such that $G$ acting on a convex polyhedron $P$ is the metric space obtained by gluing together the spherical polygons $lk(v)^*$,
spherical duals to the links of the vertices of $P$, in such a way that whenever $v_1$ and
$v_2$ share an edge, the corresponding dual edges of $lk^*(v_1)$ and $lk^*(v_2)$ are identified isometrically .}

Although, by definition, $G$ acts only on polyhedra, it is often useful to regard it as acting on the faces, edges and vertices of each polyhedron,
with their images being vertices, edges and faces of a particular cell decomposition of the corresponding metric space. This notion can be made
formal: The metric space $G(P)$ is built by gluing spherical polygons $p_i$, each corresponding to a vertex $v_i$ of $P$, and two such polygons
$p_i$, $p_j$ are glued along an edge $e'_{ij}$ if $v_i$, $v_j$ are joined by an edge $e_{ij}$. So one can define $p_i=G(v_i)$, $e'_{ij}=G(e_{ij})$,
and equivalently the point $G(f_1)$, image of a face of $P$. This picture of the generalized Gauss map as acting on the pieces rather than just the
full polyhedron will be most useful.

It can be seen that $G(P)$ is homeomorphic to $\s^2$, and that it is combinatorially Poincaré dual to $P$. If $P$ is an Euclidean polyhedron, $G(P)$
is isometric to $\s^2$. If $P$ is a hyperbolic polyhedron, a vertex of $G(P)$ is the dual to a face $f$ of $P$, $G(f)$. The edges $e_i$ surrounding
$f$ in $P$ correspond under the generalized Gauss map to edges $G(e_i)$ incident on $G(f)$. The angle at $G(f)$ between the adjacent edges $G(e_i)$
and $G(e_j)$ is equal to $\pi$ minus the face angle of $f$ between $e_i$ and $e_j$. The face $f$ is a hyperbolic polygon, so the geodesic curvature
of its boundary is, by the Gauss-Bonnet Theorem, strictly greater than $2\pi$. Therefore $G(P)$ has a \emph{cone point} at $G(f)$ with cone angle
greater than $2\pi$ and equal to the geodesic curvature of the boundary of $f$. In general, $G(P)$ has constant curvature except perhaps at a finite
number of cone points which are Poincaré dual to faces of $P$. In the next section we present an important theorem that proves that the metric on the
Gauss image of a compact convex hyperbolic polyhedron determines completely the polyhedron up to congruence. This is in contrast to Euclidean
polyhedra, whose Gauss image is always isometric to $\s^2$.

The metric on $G(P)$ can also be defined intrinsically in terms of the \emph{polar map}.

\ite{Definition}{}{The {\bf polar} $P^*$ of a convex polyhedron $P$ is the polyhedron whose boundary is the set of outward unit vectors normal
to supporting hyperplanes of $P$.}

According to the discussion on the parametrization of hyperplanes given above, the polar of an Euclidean polyhedron is always $\s^{n-1}$, while the
polars of polyhedra in $\s^n$, $\h^n$ and $\s^{n-1}_1$ are respectively polyhedra in $\s^n$, $\s^{n-1}_1$ and $\h^n$. For these three cases the polar
map defines a duality, the \emph{polar duality}, since $P^{**}=P$. Moreover, $P^*$ is a topological sphere, and it is combinatorially the Poincaré
dual of $P$, just as the Gauss image, and the metric on $P^*$ can be defined naturally as the induced metric, so that for example if $P$ is a
hyperbolic polyhedron, $P^*$ lives in $\s^2_1$ and at its smooth points the curvature is $1$. In fact the following is true \cite{Rivin1}

\ite{Proposition}{}{The intrinsic metric on the polar $P^*$ of a convex polyhedron $P$ in $\e^3$, $\s^3$ or $\h^3$ is the metric on $G(P)$.}

By the polar duality it is clear that the hyperbolic polyhedron $P$ corresponds to one and only one polyhedron $P^*$ in $\s^2_1$ under the polar map.
The question is whether the metric that arises as the Gaussian image of $P$ can be embedded isometrically as a convex polyhedron in $\s^2_1$ in two
or more (non-congruent) ways. The answer to this question is that two such non-congruent polyhedra, if they exist, must be combinatorially
inequivalent. This follows basically from the fact that the face angles of a hyperbolic polyhedron determine the polyhedron up to congruence (see
\cite{Rivin1} for details). Indeed, the metric on $P^*$ and its combinatorics determine the face angles of $P^*$, which in turn determine the planar
angles of $P$, so $P$ and $P^*$ are unique up to congruence. But the combinatorial type of $P^*$ is also determined by the metric $\cite{Rivin1}$, so
that the following holds:

\ite{Proposition\label{23}}{}{The intrinsic metric on the polar $P^*$ of a convex hyperbolic polyhedron determines $P^*$ up to congruence.}

Therefore, we shall use interchangeably the terms ``Gaussian image of $P$'', ``$G(P)$'', ``polar of $P$'', ``$P^*$'' and ``polar metric of $P$'' when
referring to the intrinsic metric on the polar of the hyperbolic polyhedron $P$.

The fact that the combinatorial type of $P^*$ is determined by its metric is an important and non-trivial result. It follows considering two
combinatorially inequivalent isometric embeddings $E_i$ of $G(P)$ as convex polyhedra, $E_1(G(P))$ and $E_2(G(P))$, with two different cell
decompositions $C_1$ and $C_2$. The cell decomposition $C$ obtained by superimposing $C_1$ and $C_2$ has the same vertices of $G(P)$ plus the
vertices where edges of $C_1$ cross edges of $C_2$. Then the polyhedra $E_1(G(P))$ and $E_2(G(P))$ can be thought of having the cell decomposition
$C$ with certain degenerate vertices and edges, so they can be thought of being combinatorially equivalent as polyhedra of some degenerate type.
These degeneracies do not spoil the argument made above for non-degenerate combinatorially equivalent polyhedra, and therefore $E_1(G(P))$ and
$E_2(G(P))$ must be congruent.

\section{Characterization of compact convex hyperbolic polyhedra}

\setcounter{sec}{1} \label{compact}

In this section we present a theorem that gives a characterization of compact convex hyperbolic polyhedra in terms of the metric on the polar. A
detailed exposition together with a complete proof can be found in \cite{Rivin1}. We begin with a definition.

\ite{Definition \label{31}}{}{A metric space $Q$ is called {\bf admissible} if:
\begin{enumerate}
\item $Q$ is homeomorphic to the sphere $\s^2$.
\item $Q$ is piecewise spherical, with constant curvature $1$ away from a set of cone points.
\item All cone angles are $>2\pi$.
\item Every closed geodesic on $Q$ has length $>2\pi$.
\end{enumerate}}

Admissible metrics are precisely the images under the generalized Gauss map of compact convex hyperbolic polyhedra. Moreover, such a polyhedron is
determined up to congruence by this metric. This is the content of Rivin's theorem for the characterization of compact convex polyhedra
\cite{Rivin1}.

\ite{Theorem \label{32}}{(Rivin)}{A metric space $Q$ is the polar of a unique (up to congruence) compact convex polyhedron in $\h^3$ if and only if
$Q$ is admissible.}

This theorem provides a one-to-one correspondence between admissible metric spaces and compact convex polyhedra. More specifically this can be stated
in following way. Let $\mathcal{P}_n$ be the space of compact convex hyperbolic polyhedra with $n$ numbered faces, modulo congruence. Let
$\mathcal{M}_n$ be the space of admissible metrics on $\s^2$ with $n$ numbered cone points, modulo isometry. The Gauss map defines a map
$\phi:\mathcal{P}_n\to\mathcal{M}_n$ which takes a hyperbolic polyhedron to the intrinsic metric on its polar. Then Rivin's theorem can be stated
as:\\

\n{\bf Theorem.} \emph{The map $\phi$ is a bijection.}\\

In order to prove this theorem one must show that: 1.- $\phi(\mathcal{P}_n)$ is in $\mathcal{M}_n$. This proves that the polar of a compact
hyperbolic polyhedron is an admissible metric. 2.- $\phi$ is one-to-one. This is equivalent to the fact that non-congruent polyhedra have
non-isometric polar metrics. And 3.- $\phi$ is surjective. This proves that all admissible metrics arise as the Gaussian image of a compact convex
polyhedron in $\h^3$.

For the proof of 1 we first present some facts concerning polygonal curves in hyperbolic and de Sitter spaces.

\ite{Lemma \label{33}}{}{Let $\gamma$ be a closed geodesic on the polar $P^*\subset \s^2_1$ of a compact convex hyperbolic polyhedron. Then the polar
of $\gamma$ in $\h^3$ is the set of edges of an infinite polyhedral cylinder $C$ in $\h^3$. Moreover, the length of $\gamma$ is equal to the sum of
the exterior dihedral angles of $C$.}

\proof The polygonal curve $\gamma$ can be split into geodesic segments $\gamma_i$ in $\s^2_1$, joined at edges or vertices of $P^*$. The polar of
each $\gamma_i$ is a geodesic $\Gamma_i$ in $\h^3$. If $\gamma_i$ and $\gamma_j$ meet at a point $p$ on a edge or a vertex of $P^*$, then by the
polar duality $\Gamma_i$ and $\Gamma_j$ contain two sides of the face $f$ polar to $p$, so $\Gamma_i$ and $\Gamma_j$ are coplanar. Moreover,
$\Gamma_i$ and $\Gamma_j$ cannot meet in $\h^3$ or even on the sphere at infinity, since the meeting point would correspond, under the polar map, to
a face of $P^*$ containing both $\gamma_i$ and $\gamma_j$. This could only happen if $\gamma_i$ and $\gamma_j$ are adjacent sides of a face in $P^*$,
but then they would subtend an angle less than $\pi$ (the face angle at $p$), which contradicts the fact that $\gamma$ is a geodesic. Therefore
$\Gamma_i$ and $\Gamma_j$ are coplanar and hyperparallel. But the $\gamma_i$ are joined cyclically, so the $\Gamma_i$ form the edges of an infinite
polyhedral cylinder $C$. Let the two adjacent faces (infinite geodesic strips) to $\Gamma_i$ in $C$ be $f_1$ and $f_2$. Under polar duality $f_1$ is
mapped to the initial point of $\gamma_i$ in $P^*$ and $f_2$ to the end point, and $\Gamma_i$ is mapped to $\gamma_i$ in such a way that the exterior
dihedral angle between $f_1$ and $f_2$ is equal to the length of $\gamma_i$. So the length of $\gamma$ is equal to the sum of the exterior dihedral
angles of $C$.\qed

\ite{Lemma \label{34}}{}{Let $H_1$ and $H_2$ be two geodesic half planes meeting at an (exterior) dihedral angle $\alpha$. Let $\gamma$ be a curve on
$H_1\cup H_2$ which is a geodesic in the intrinsic metric on $H_1\cup H_2$, and intersecting $H_1\cap H_2$ at a (unique) point $p$. Then the turning
angle of $\gamma$ at $p$ is $\tau_p(\gamma)\le \alpha$.}

\proof Consider a unit sphere centered at $p$. The intersection of $H_1\cup H_2$ with the sphere is a lune with internal angle $\pi-\alpha$. Consider
one vertex of the lune, $q$, and the two intersection points $p_1$ and $p_2$ of $\gamma$ with the sphere. Clearly $p_1$ and $p_2$ are on opposite
edges of the lune. Then $qp_1p_2$ defines a spherical triangle with an angle equal to $\pi-\alpha$, and by Snell's law the adjacent sides are of size
$\beta$ and $\pi-\beta$, where $\beta$ is the angle between $\gamma$ and $H_1\cap H_2$. Then the turning of $\gamma$ at $p$ is $\pi$ minus the length
$\ell$ of the side $p_1p_2$ of this spherical triangle, which is, by the spherical law of cosines,
\eqa{\cos\ell&=&\cos\beta\cos(\pi-\beta)+\sin\beta\sin(\pi-\beta)\cos(\pi-\alpha)\nonumber\\[2mm]
&=&-\cos^2\beta+\sin^2\beta\cos(\pi-\alpha)\le \cos(\pi-\alpha).\nonumber}
Therefore $\ell\ge(\pi-\alpha)$, so that $\tau_p(\gamma)=(\pi-\ell)\le\alpha$. \qed

\ite{Theorem \label{35}}{(Hyperbolic Frenchel Theorem)}{The total turning of a closed polygonal curve $\gamma$ in $\h^3$ not contained in a geodesic
is $\tau(\gamma)>2\pi$.}

\proof The fact follows from considering the set of triangles $T_i=p_1p_ip_{i+1}$, with $2\le i\le k-1$, that form an immersed disk with boundary
$\gamma$. Then the sum of the angles $\alpha_i^j$ of the $T_i$ meeting at a vertex $p_j$ in $\gamma$ is $\sum_i\alpha_i^j\ge\pi-\tau_j(\gamma)$,
where $\tau_j(\gamma)$ is the turning of $\gamma$ at $p_j$ (equality holds when the $T_i$ meeting at $p_j$ are coplanar, in particular for $p_2$ and
$p_{k}$, and the inequality follows from the spherical triangle inequality). The total turning of $\gamma$ is then
$\tau(\gamma)=\tau_1(\gamma)+\sum_{j=2}^k\tau_j(\gamma)\ge\tau_1(\gamma)+(k-1)\pi-\sum_{j=2}^k\sum_i\alpha_i^j$. The turning at $p_1$,
$\tau_1(\gamma)$ can be bounded also with the spherical triangle inequality taking the triangles in pairs, which gives
$\tau_1(\gamma)\ge\pi-\sum_{i=2}^{k-1}\alpha_i^1$. Therefore, $\tau(\gamma)\ge k\pi-\sum({\rm all\ angles\ of\ all\ }T_i)>k\pi-(k-2)\pi=2\pi$, since
the angles of a hyperbolic triangle sum up lo less than $\pi$. So $\tau(\gamma)>2\pi$. \qed

This brings us to the first part of the proof of Rivin's theorem:

\ite{Proposition \label{36}}{}{The Gaussian image of a compact convex hyperbolic polyhedron $P$ is an admissible metric.}

\proof First, it is clear from the discussion in Section~\ref{defs} that $G(P)$ is homeomorphic to $\s^2$, with constant curvature $1$ away from a
set of cone points polar to the faces of $P$, with cone angles $>2\pi$. It remains to prove that all closed geodesics in $P^*$ are longer than
$2\pi$. Let $\gamma$ be such a geodesic. By Lemma~\ref{compact}.3, $\ell(\gamma)=\sum_i \alpha_i$ where $\alpha_i$ are the exterior dihedral angles
of the corresponding cylinder $C$. Let $\rho$ be a meridian curve on $C$. By Lemma~\ref{compact}.4, $\tau(\rho)\le \sum_i\alpha_i$, and by
Theorem~\ref{compact}.5, $\tau(\rho)>2\pi$. Therefore $\ell(\gamma)=\sum_i \alpha_i\ge\tau(\rho)>2\pi$, so $\ell(\gamma)>2\pi$.\qed

The second part of Rivin's theorem (uniqueness), is the content of the following theorem, which follows from Proposition~\ref{defs}.5 and the duality
of the polar map.

\ite{Theorem}{}{The intrinsic metric on the polar of a compact convex hyperbolic polyhedron determines the polyhedron up to congruence.}

Finally one must prove that all the admissible metrics in $\mathcal{M}_n$ are the image under $\phi$ of a compact convex hyperbolic polyhedron in
$\mathcal{P}_n$. Here we sketch the proof by Rivin and Hodgson \cite{Rivin1}.

First, the space of metrics on a compact space $M$ (that is, the set of metric spaces $M_i$ homeomorphic to $M$) can be endowed with a distance (the
\emph{Lipschitz distance}). Basically, the Lipschitz distance between $M_1$ and $M_2$ is $d_L(M_1,M_2)={\rm inf}\ \mathcal{D}(f)$, where $f$ ranges
over all homeomorphisms from $M_1$ to $M_2$, and $\mathcal{D}(f)$ is a measure of the largest distortion in distance that the map $f:M_1\to M_2$
produces. In particular, if two metric spaces are isometric then there is an isometry (with zero distortion) which can be chosen as $f$, so the
Lipschitz distance is zero. Then the space $\mathcal{M}_n$ can be given a metric, the \emph{Lipschitz metric}, in much the same way but $f$ ranging
over homeomorphisms that take cone points to cone points preserving the numberings.

Then, it can be proven that any two admissible metrics $g_0$ and $g_1$ with a finite number of cone points can be connected by a path $g_t$,
$t\in[0,1]$, of admissible metrics, which is continuous in the Lipschitz topology. Moreover, all the metrics $g_t$ with $t\in(0,1)$ have the same
number of cone points $N$, so a continuous deformation of a metric can only change the number of cone points in an attained limit. Also, if $g_0$ is
the polar of a compact convex polyhedron (that is, belongs to $\phi(\mathcal{P}_n)$ for some $n$), then this is also true for all $g_t$ with $t$
sufficiently close to $0$.

The space $\widetilde{\mathcal{P}}_n$ of all convex hyperbolic polyhedra with $n$ numbered faces endowed with the Haussdorf topology is a
$3n$-dimensional manifold, since it can be (openly) embedded in $(\s^2_1)^n$ (polars of faces). The space
$\mathcal{P}_n=\widetilde{\mathcal{P}}_n/{\rm isom}(\h^3)$ is therefore an open manifold of dimension $3n-6$, since the group of isometries acts
freely and properly on $\widetilde{\mathcal{P}}_n$. Also, the space $\mathcal{M}_n$ is a manifold of dimension $3n-6$. This can be seen by
considering the space $\mathcal{S}_n$ of piecewise spherical metrics on $\s^2$ with $n$ labeled points (not necessarily admissible) up to isometry.
Choosing a suitable triangulation, and noting that by Euler's formula any triangulation with $n$ vertices has $3n-6$ edges, $\mathcal{S}_n$ can be
embedded into $\R^{3n-6}$ (by considering the lengths of the $3n-6$ edges). The conclusion follows from the fact that the conditions of admissibility
in $\mathcal{S}_n$ are open conditions. Now, since $\mathcal{P}_n$ and $\mathcal{M}_n$ are open manifolds and $\phi:\mathcal{P}_n\to \mathcal{M}_n$
is a continuous injective map, then by the theorem of invariance of domain it is an open map. This proves that if $g(t_0)\in\phi(\mathcal{P}_n)$ for
some $t_0$, then $g(t)\in\phi(\mathcal{P}_n)$ for $t$ in a neighborhood of $t_0$.

Now consider a sequence $\{P_i\}$ of compact hyperbolic polyhedra, and consider the corresponding sequence of admissible metrics $\{g_i\}$, and
assume that the sequence $\{g_i\}$ converges in the Lipschitz topology to a metric $g_\infty$. Then it can be proven that if $g_\infty$ is admissible
($g_\infty\in\mathcal{M}_n$), then $g_\infty\in\phi(\mathcal{P}_n)$. Indeed, if the $P_i$ have bounded diameters, then there is a subsequence that
converges to a compact convex limit $P_\infty$ with polar metric $g_\infty$. If $g_\infty\notin\phi(\mathcal{P}_n)$ then it is because $\{P_i\}$
degenerates in the limit to a compact polygon, segment or to a point. In those cases $g_\infty$ is a \emph{suspension}, with two antipodal cone
points with cone angle $>2\pi$, and any geodesic passing through both cone points has length $2\pi$, so $g_\infty$ is not admissible. If the $P_i$ do
not have bounded diameters then there is a subsequence with ${\rm diam}(P_i)\to\infty$. Every compact hyperbolic polyhedron $P$ with $n$ vertices and
of sufficiently large diameter (how large depends on $n$) contains what is called a \emph{long thin tube}. This is an arbitrarily thin polyhedral
cylinder $C$ such that there is a hyperbolic plane arbitrarily far from the vertices of $P$ and intersecting $C$ in a way arbitrarily close to
perpendicular. The edges of $C$ are also arbitrarily close to parallel between them. (For example, one may look at the limiting case of a polyhedron
$P'$ with an ideal vertex $v$, in the Poincaré model. A small enough sphere orthogonal to the sphere at infinity and enclosing $v$ is a hyperbolic
plane that intersects $P'$ across a long thin tube.) Therefore, the polar of $C$ is a piecewise geodesic curve $\gamma$ in $P^*$ with total turning
arbitrarily close to 0, and with length arbitrarily close to $2\pi$ (because the intersection of the hyperplane with $C$ is a polygon arbitrarily
close to Euclidean). In the limit this turns into a geodesic of length $2\pi$ in $g_\infty$, so $g_\infty$ is not admissible.

With these considerations, one can prove the following.

\ite{Proposition}{}{Every admissible metric is the metric on the polar of a compact convex hyperbolic polyhedron.}

Indeed, choose an admissible metric $g'\in\mathcal{M}_n$. According to the previous discussion, there is a continuous path $g:[0,1)\to\mathcal{M}_N$
such that $g(t)\in\phi(\mathcal{P}_N)$  for all $t$ close to $0$, and $\lim_{t\to 1}g(t)=g'$. Then by the same arguments, $g(t)$ is in fact inside
$\phi(\mathcal{P}_N)$ for all $t\in[0,1)$. Finally, since $g'$ is admissible and it is the limit of a sequence of admissible metrics $g(t)$, then
$g'\in\phi(\mathcal{P}_n)$, so it is the metric on the polar of a compact convex polyhedron with $n$ faces.

\section{Characterization of ideal polyhedra}

\label{ideal} \setcounter{sec}{1}

Ideal polyhedra can be characterized in an equivalent way as compact polyhedra, that is, in terms of the polar metric. It is clear from Rivin's
theorem that the polar of an ideal polyhedra cannot be an admissible metric. This means that it is necessary to enlarge the set of metrics on $\s^2$
to include those which arise as the Gaussian image of ideal polyhedra. We will see that they form some sort of ``boundary'' of the set of admissible
metrics.

\ite{Definition}{}{A metric space $Q$ is called {\bf ideally admissible} if:
\begin{enumerate}
\item $Q$ is homeomorphic to the sphere $\s^2$.
\item $Q$ has the structure of a cell complex combinatorially equivalent to a convex polyhedron.
\item Each 2-cell is isometric to a hemisphere of the standard unit sphere $\s^2$.
\item Any two cells $H$ and $K$ are either disjoint or they are identified along an equatorial arc of length $\ 0<l(H,K)<\pi$.
\item The length of any simple cycle in the 1-skeleton of $Q$ is strictly greater than $2\pi$ unless it forms the boundary of a 2-cell.
\end{enumerate}}

In the same way that the Gauss image of a compact convex polyhedron is the topological sphere made up gluing convex spherical polygons, corresponding
to the (finite) vertices of the polyhedron, the polar of an ideal polyhedron is made up of a gluing of hemispheres. The main theorem concerning the
characterization of ideal polyhedra is due to Rivin \cite{Rivin2,Rivin4}.

\ite{Theorem}{(Rivin)}{A metric space $Q$ is the polar of an ideal convex hyperbolic polyhedron $P$ if and only if $Q$ is ideally admissible.}

\begin{figure}
\begin{center}
\includegraphics[height=5cm,width=14cm]{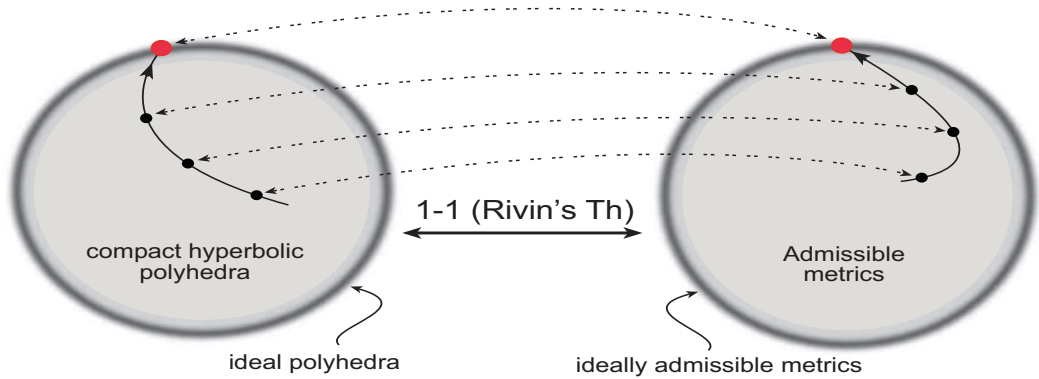}
\end{center}
\caption{Ideally admissible metrics arise as limits of sequences of admissible metrics. Rivin's Theorem then implies that ideal polyhedra arise as
limits of sequences of compact convex polyhedra.} \label{f4.1}
\end{figure}

The main idea underlying Rivin's proof of this theorem is that ideally admissible metrics arise as the limit of certain sequences of admissible
metrics. According to Rivin's theorem of compact polyhedra, these sequences are in one-to-one correspondence with sequences of compact polyhedra (see
Fig.~\ref{f4.1}). Therefore, in order to prove the theorem one must show that: 1.-the limit of such sequences of compact polyhedra are ideal
polyhedra, 2.-this limit commutes with the Gauss map, and 3.-the Gaussian image of an ideal polyhedron is an ideally admissible metric.

\ite{Definition}{}{Let $Q$ be an ideally admissible metric space, and let $t>0$. We define $Q^t$ (the t-expansion of $Q$) as the metric space
constructed in the following way:
\begin{enumerate}
\item Divide each 2-cell $H$ (a spherical hemisphere) of $Q$ into triangles with a vertex in the center of $H$, adjacent sides being radii of
length $\pi/2$, and opposite sides the edges of the 1-skeleton of $Q$ in $\partial H$. $Q$ can be thought of as the particular gluing $\mathcal{G}$
of these triangles.
\item Substitute each triangle by a spherical triangle exactly the same, except that the length of the side on $\partial H$ is
multiplied by $(1+t)$.
\item Glue these triangles back again following the pattern of $\mathcal{G}$ to form the metric space $Q^t$.
\end{enumerate}}

\begin{figure}
\begin{center}
\psfrag{Q}{$Q$} \psfrag{T}{$Q^t$} \psfrag{l1}{$\ell_1$} \psfrag{l2}{$\ell_2$} \psfrag{l3}{$\ell_3$} \psfrag{t1}{\hspace{-1.1cm}$(1+t)\ell_1$}
\psfrag{t2}{\hspace{-0.7cm}$(1+t)\ell_2$} \psfrag{t3}{\hspace{-0.5cm}$(1+t)\ell_3$}
\includegraphics[height=5cm,width=14cm]{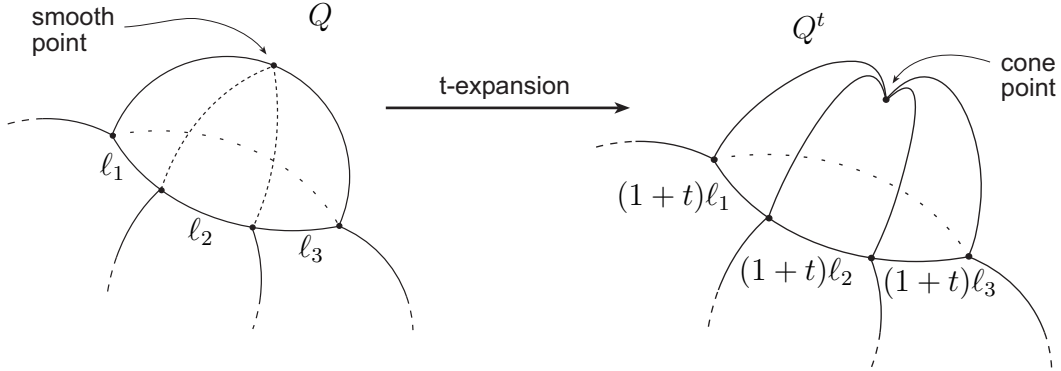}
\end{center}
\caption{t-expansion of an admissible metric. Each $2$-cell (a spherical hemisphere) is substituted by a hemisphere of constant curvature $1$ with an
equator of length $2(1+t)\pi$. Therefore a cone point arises at the north pole (analogous to a \emph{suspension}) with cone angle $2(1+t)\pi$.}
\label{f4.2}
\end{figure}

The construction of $Q^t$ from an ideally admissible metric space $Q$ is exemplified in Fig.~\ref{f4.2}. Let $e_1(Q)$ denote the longest edge of the
1-skeleton of $Q$. The following theorem is the starting point of the program outlined above.

\ite{Theorem}{}{Let $0<t<\pi/e_1(Q)-1$. Then $Q^t$ is admissible.}

\proof First, by construction, $Q^t\simeq Q \simeq \s^2$ and $Q^t$ is a spherical cone manifold. Second, the cone angles of $Q^t$ are strictly
greater than $2\pi$. In fact, there are two types of cone points: those corresponding to vertices of $Q$ (\emph{ordinary} vertices), with cone angles
of at least $3\pi$, and those arising at the center of each hemisphere of $Q$ through the $t$-expansion (\emph{special} vertices) as sketched in
Fig.~\ref{f4.2}. By construction special vertices have cone angles equal to $2(1+t)\pi$. It remains to see that closed geodesics on $Q^t$ are
strictly longer than $2\pi$. We shall be using the term \emph{face} of $Q^t$ for the star of a special vertex, that is, a hemisphere of $Q$ after the
$t$-expansion.

Let $\gamma_H$ be a nonempty connected component of the intersection of a closed geodesic with the interior of a face $H$ of $Q^t$. If $\gamma_H$
does not pass through the cone point of $H$, then together with the boundary of $H$ bounds a spherical lune at one of its sides, so
length$(\gamma_H)=\pi$. If $\gamma_H$ passes through the cone point, then it is composed of two spherical radii, so length$(\gamma_H)=\pi$. This
shows that if a closed geodesic intersects the interior of three faces of $Q^t$ then its length is at least $3\pi$. Suppose now that a closed
geodesic $\gamma$ is composed of two of such geodesic arcs, then the edge $E_t$ joining the two faces must be at least $\pi$. This edge is of length
$\ell(E_t)=(1+t)\ell(E)$, where $E$ is the edge before the $t$-expansion, so $\ell(E_t)\le(1+t)e_1(Q)$. By hypothesis $t<\pi/e1(Q)-1$, so
$\ell(E_t)\le\pi$, which proves that $\gamma$ cannot be composed of two geodesic arcs intersecting the interior of two faces.

Now consider a closed geodesic $\gamma$ that does not intersect the interior of any face of $Q^t$. Then $\gamma$ lives in the $t$-expansion of the
1-skeleton of $Q$. By definition of $Q$ the length of this geodesic is $\ell(\gamma)\ge 2(1+t)\pi>2\pi$.

Finally, consider the remaining case of a closed geodesic $\gamma$ that intersects exactly one face $H$ of $Q^t$. Let $p$ and $q$ be the intersection
of $\gamma$ with $\partial H$. Let $\rho_1$ and $\rho_2$ be the two segments from $p$ to $q$ along $\partial H$, and let $\gamma_1$ and $\gamma_2$ be
$\gamma$ with $\gamma_H$ replaced with $\rho_1$ and $\rho_2$ respectively ($\gamma_i$ are geodesics). If both $\gamma_1$ and $\gamma_2$ bound a face
of $Q^t$ then $Q$ must be composed of just 3 cells, and then it cannot be ideally admissible (all the gluings cannot be shorter than $\pi$). So at
least one of them, say $\gamma_1$, does not bound a face of $Q^t$, and therefore $\ell(\gamma_1)>2(1+t)\pi$. Moreover, from
$\ell(\rho_1)+\ell(\rho_2)=2(1+t)\pi$ and $\ell(\rho_i)\ge\pi$ it follows that $\ell(\rho_1)\le \pi+2t\pi$, so
\eq{\ell(\gamma)=\ell(\gamma_1)-\ell(\rho_1)+\ell(\gamma_H)=\ell(\gamma_1)-\ell(\rho_1)+\pi>(2(1+t)\pi)-(\pi+2t\pi)+\pi=2\pi}
and therefore $\ell(\gamma)>2\pi$. \qed

In order outline the proof of parts 1 and 2, we consider certain sequences $\{t_j\}_{j=1}^\infty$ with $\lim t_i=0$, and the corresponding sequences of polyhedra $P^t$ with polars $Q^t$. First, it can be seen that the set of edges of the cell decomposition of $Q^t$ as constructed in Definition 4.3 is the same as that of the cell decomposition of the Gauss image of $P^t$ as constructed in Definition 2.2, so the following is true about $P^t$. Each face of $P^t$ corresponds to a vertex of $Q^t$, and it is called \emph{ordinary} or \emph{special} depending on whether the corresponding vertex in $Q^t$ is an ordinary or a special vertex. Then, special faces of $P^t$ are surrounded completely by ordinary faces and the dihedral angles between every special face and the ordinary faces are always $\pi/2$. On the other hand, ordinary faces are surrounded by an alternating sequence of ordinary and special faces, and the dihedral angle between two ordinary faces $F_i$ and $F_j$ is $(1+t)\ell(F^*_i,F^*_j)$, where $\ell(F^*_i,F^*_j)$ is the length of the edge joining the (polar) vertices $F^*_i$ and $F^*_j$ in $Q$.

\begin{figure}
\begin{center}
\includegraphics[height=5cm,width=8cm]{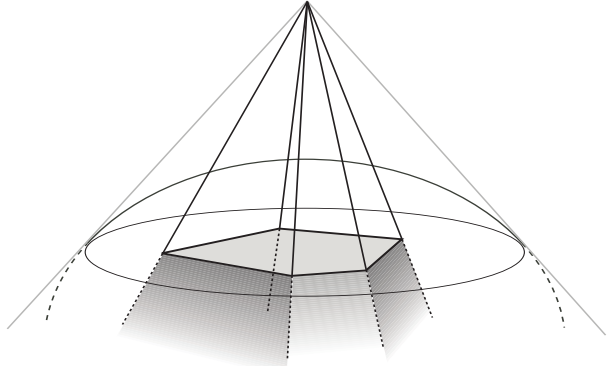}
\end{center}
\caption{A special face of the compact convex polyhedron $P^t$ has dihedral angles of $\pi/2$ with all its surrounding (ordinary) faces. Therefore it
can be represented by a hyperinfinite vertex.} \label{f4.3}
\end{figure}

Since special faces are orthogonal to all their (ordinary) surrounding faces, these ordinary faces might be extended geodesically to meet at a
hyperinfinite vertex beyond the sphere at infinity. Therefore special faces can be represented by hyperinfinite vertices, and the polyhedra $P^t$ can
be thought of polyhedra with all their vertices hyperinfinite (see Fig.~\ref{f4.3}). The cone angle at a special vertex in $Q^t$ is equal to the sum
of the exterior face angles of the corresponding special face of $P^t$. Since $Q^t\to Q$, the cone angle goes to $2\pi$ in the limit, so the special
faces turn Euclidean in the limit $t\to 0$, and hence of zero area.

It can be shown by arguments on elementary hyperbolic plane geometry \cite{Rivin2}, that for $t$ small enough there exists a ball $B_R$ of radius $R<\infty$ in $\h^3$ which intersects every ordinary face of $P^t$  and every edge between pairs of such faces. Moreover, the set of planes in $\h^3$ intersecting such ball is always compact. Therefore, we may consider the planes $\Pi_i^t$ containing the ordinary faces of $P^t$, and a sequence $\{t_j\}_{j=1}^\infty$  with $\lim t_k=0$ such that every sequence $\{\Pi_i^{t_j}\}$ converges to a plane $\Pi_i$ that intersects $B_R$. Also, by continuity, the dihedral angle formed by $\Pi_i$ and $\Pi_j$ is exactly $\ell(F^*_i,F^*_j)$.

Now we consider the intersection of the oriented half-spaces defined by these $\Pi_i$, which defines a generalized polyhedron $\overline{P}$. Again, by continuity, $\overline{P}$ cannot  have finite vertices, since the position of such vertices is the limit of the hyperinfinite vertices of $P^t$. The claim is that the dual $\overline{P}^*$ is isometric to $Q$, and $\overline{P}$ it is an ideal polyhedron. Indeed, each cycle in the one-skeleton of $\overline{P}^*$ has length $2\pi$, and such cycles are geodesics in $\overline{P}$. Such cycle corresponds to an ideal vertex of $\overline{P}$, and therefore it bounds a round hemisphere.

Now we prove part 3, i.e. that the Gaussian image of an ideal polyhedron is an ideally admissible metric.

\ite{Theorem}{}{The metric space $P^*$ arising as the Gaussian image of a convex ideal polyhedron $P$ is ideally admissible.}

\proof $P^*$ is the metric space homeomorphic to $\s^2$ obtained from the gluing of the spherical polygons polar to the vertices of $P$. In order to
prove that $P^*$ is ideally admissible we must prove that: {(a) The} spherical polygons dual to vertices of $P$ are spherical hemispheres. \mbox{(b)
The} gluing edges have lengths in the interval $(0,\pi)$. (c) The simple cycles in $P^*$ not bounding a face are strictly longer than $2\pi$.

(a) Let $v$ be an (ideal) vertex of $P$. The perimeter or $v^*$ is equal to the sum of the exterior dihedral angles at the edges incident on $v$ in
$P$, which is $2\pi$. The classical way to see this is to consider a small horosphere tangent to the sphere at infinity at $v$. This horosphere has
Euclidean intrinsic metric, so the intersection of $P$ with the horosphere is an Euclidean convex polygon with exterior angles equal to the exterior
dihedral angles of $P$ at $v$, and summing up to $2\pi$ (like any convex Euclidean polygon). Then, since the perimeter of the convex spherical
polygon $v^*$ is $2\pi$, $v^*$ is a spherical hemisphere.

(b) Since $P$ is convex, its dihedral angles are contained in the interval $(0,\pi)$, so the length of the edges of $P^*$ are also in the same
interval.

(c) This proof is very similar to the one for compact polyhedra of Section~\ref{compact}. A simple cycle in $P^*$ not bounding a face correspond to a
chain $F$ of faces of $P$ joined by edges $e_i$ all of which do not share a common ideal vertex. $F$ can be extended geodesically beyond its
boundaries to form a complete hyperbolic surface $\tilde{F}$ immersed in $\h^3$. Now, $\tilde{F}$ is topologically an infinite cylinder with both
edges of infinite volume. Let $\gamma$ be a closed geodesic on $\tilde{F}$ homotopic to the meridian curve, which is embedded in $\h^3$ as a
polygonal curve with turnings at the edges $e_i$. By the Hyperbolic Frenchel's Theorem the total turning of $\gamma$ is greater than $2\pi$, and by
Lemma~\ref{compact}.4 the sum of dihedral angles at $e_i$ is greater or equal to the total turning of $\gamma$. So the length of the corresponding
cycle in $P^*$ is larger than $2\pi$. \qed

In the projective model of $\h^3$, ideal polyhedra are represented by Euclidean polyhedra with all the vertices on the sphere at infinity. Therefore
there is a one-to-one correspondence between hyperbolic ideal polyhedra and Euclidean polyhedra inscribed in the sphere. The characterization of
ideal hyperbolic polyhedra then solves the problem (posed by J.~Steiner almost two centuries ago) of giving a combinatorial characterization of
Euclidean polyhedra inscribed in the sphere. Moreover, given a polyhedron $P$, there is a polynomial-time algorithm that decides whether $P$ can be realized as an ideal hyperbolic polyhedron, and therefore as an Euclidean polyhedron inscribed in the sphere \cite{Rivin3}.  An example of a class of non-inscribable polyhedra is given in \cite{Rivin4}, in terms of the \emph{stellation} of certain polyhedra. The \emph{stellation} of a polyhedron $P$ is defined as the polyhedron obtained from $P$ by adding a vertex $v$ for each face $f$ of $P$, and replacing every such face by the union of the triangles with apex $v$ and base the edges of $f$. 

\ite{Theorem}{}{Let $P$ be a polyhedron and let $V(P)$ and $F(P)$ denote the number of vertices and faces of $P$. Then, if $V(P)\le F(P)$ the stellation of $P$ cannot be inscribed in the sphere.}

\proof Consider a polyhedron $P$ whose stellation $S$ can be inscribed in the sphere. We consider the ideal polyhedron $S_I$ associated with $S$. Its vertices can be divided in two classes: the set $V_P$ of vertices original of $P$ and the set $V_S$ of vertices arising in the process of stellation.  Consider a vertex $v$ in $S_I$, and all the edges incident to $v$. The sum of the external dihedral angles at all these edges is denoted by $d(v)$. Since all the vertices are ideal, we have that $d(v)=2\pi$ for all $v$. Therefore, $\sum_{v\in V_P}d(v)=2\pi V(P)$ and $\sum_{v\in V_S}d(v)=2\pi F(P)$. But all the edges incident to vertices in $V_S$ are also incident to vertices in $V_P$, so necessarily $\sum_{v\in V_P}d(v)>\sum_{v\in V_S}d(v)$. This proves that if $S$ is inscribable then it must hold that $V(P)>F(P)$. \qed
 
In fact it is very easy to find such polyhedra: choose a random abstract polyhedron $P$; if $V(P)\le F(P)$ then $P$ is it, and if $V(P)> F(P)$ then take the Poincaré dual. For example a cube is not such a polyhedron, since it has $8$ vertices and $6$ faces, but its Poincaré dual, the octahedron, has a non-inscribable stellation. In general any polyhedron with triangular faces is of this kind. This can be seen in the following way: If $P$ has triangular faces then $E(P)=3F(P)/2$, so by Euler's formula $F(P)=2V(P)-4$, which for $V(P)>3$ means that $F(P)\ge V(P)$.

\section{Adreev's theorem from Rivin's characterization}

\setcounter{sec}{1}

\begin{figure}
\begin{center}
\includegraphics[height=5cm,width=14cm]{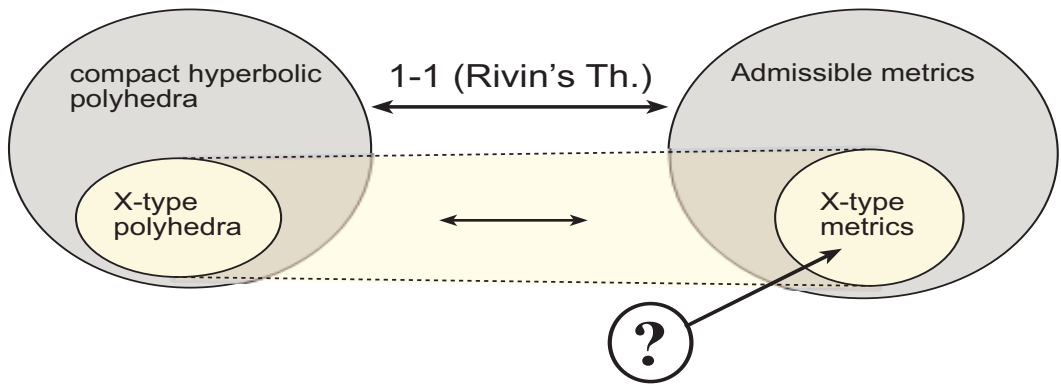}
\end{center}
\caption{Rivin's Theorem says that for each compact convex hyperbolic polyhedron there is one and only one admissible metric. But what kind of
admissible metrics are associated with a specific type $X$ of polyhedra?} \label{f5.1}
\end{figure}

Rivin's theorem provides a characterization of convex hyperbolic polyhedra in terms of the polar, which is an admissible metric on the 2-sphere. In
particular, it establishes a one-to-one correspondence between polyhedra in $\h^3$  and admissible metrics on $\s^2$. It does not, however, say much
about how specific subsets of the space of hyperbolic polyhedra are mapped under the polar map. The characterization of polyhedra with such specific
properties is given, according to Rivin's theorem, in terms of a particular subset of the space of admissible metrics, but the identification of such
subset must be built on top of it. The situation is sketched in Fig.~\ref{f5.1}.

Nevertheless, Rivin's theorem does simplify considerably the proofs of characterization theorems for specific polyhedral types. We have already seen
an example of this procedure in Section~\ref{ideal}, concerning the characterization of ideal polyhedra. In that case, ideal polyhedra were obtained
as the dual of a polar metric $Q$, constructed as the limit of a sequence of admissible metrics $Q^t$ dual to polyhedra with finite vertices, where
Rivin's theorem is at work.

A general procedure to characterize specific types of hyperbolic polyhedra that benefits from the power of Rivin's theorem is the following. Suppose
we want to prove a theorem that characterizes hyperbolic polyhedra of type $X$. Then this theorem would be stated as:\\

\emph{`` {\bf Theorem X:} An abstract polyhedron $P$ is realizable as a convex polyhedron of type $X$ in $\h^3$ if and only if the conditions $Y$ are
satisfied.''}\\

Now, from Rivin's theorem we know that $P$ is realizable as a convex polyhedron of type $X$ if and only if its polar $P^*$ defines an admissible
metric on $\s^2$ of some type $X$. Therefore, the steps to follow in order to prove Theorem X are:

\begin{enumerate}
\item Translate conditions $Y$ into conditions on the polar $P^*$ (conditions $X$ on the set of metrics on $\s^2$) defining what we may call
(conveniently) a set of $X$-metrics on $\s^2$.

\item Prove that the polar of a polyhedron of type $X$ is an $X$-metric.

\item Prove that the set of $X$-metrics is contained in the set of admissible metrics.
\end{enumerate}

Step~1 formulates the problem in the context of polars so we can use the polar duality. Then the conditions $Y$ are proven necessary in step~2 and
sufficient in step~3.

A famous example of a ``theorem X'' is given by Andreev's theorem \cite{Andreev} for the characterization of convex hyperbolic polyhedra with
dihedral angles $\le \pi/2$. The procedure outlined above was followed by Hodgson \cite{Hodgson} to derive Andreev's theorem from Rivin's
characterization (see Fig.~\ref{f5.2}). Before stating Andreev's theorem we need a definition:

\ite{Definition}{}{Let $P$ be a convex polyhedron. A {\bf k-prismatic element} is a circular sequence of $k$ edge-adjacent faces of $P$ such that no
three of these faces have a common point.}

Andreev's theorem is then stated as follows:

\begin{figure}
\begin{center}
\includegraphics[height=5cm,width=14cm]{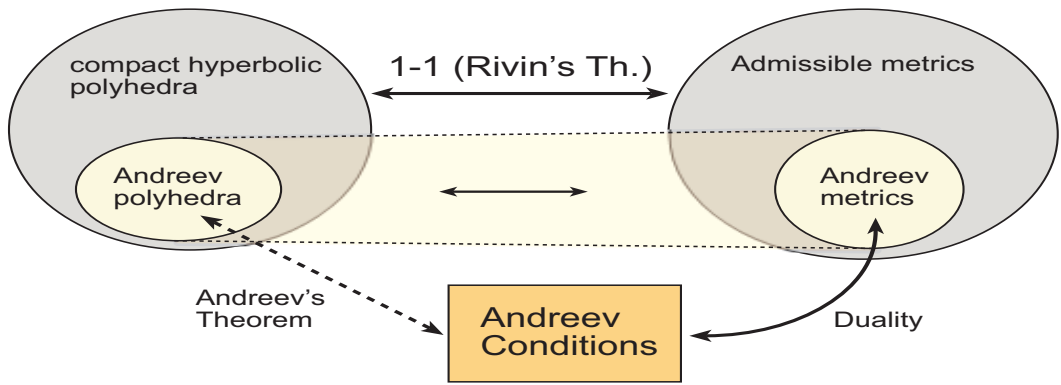}
\end{center}
\caption{Sketch of how Rivin's Theorem can be used to characterize Andreev polyhedra: One sets a dictionary between Andreev's conditions and
conditions on the space of admissible metrics, and Rivin's Theorem does the rest of the work.} \label{f5.2}
\end{figure}

\begin{figure}
\begin{center}
\psfrag{p}{$\pi/2$}
\includegraphics[height=4cm]{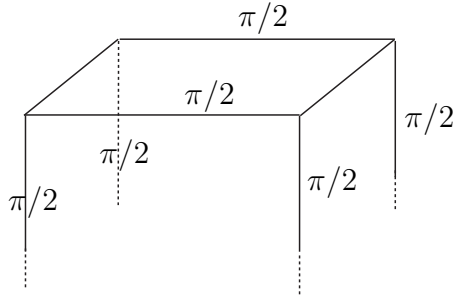}\hspace{1cm}
\end{center}
\caption{A configuration that cannot be part of a compact hyperbolic polyhedron with all dihedral angles less of equal to $\pi/2$.} \label{f5.3}
\end{figure}

\ite{Theorem}{(Andreev)}{Let $P$ be a compact convex polyhedron in $\h^3$ with faces $F_i$ and dihedral angles $\alpha_{ij}\le \pi/2$ between faces
$F_i$ and $F_j$. Then $P$ has trivalent vertices, and
\begin{enumerate}
\item \ $0<\alpha_{ij}\le \pi/2$.
\item \ If $F_i\cap F_j\cap F_k$ is a vertex then $\alpha_{ij}+\alpha_{jk}+\alpha_{ki}>\pi$.
\item \ If $F_i$, $F_j$, $F_k$ form a 3-prismatic element, then $\alpha_{ij}+\alpha_{jk}+\alpha_{ki}<\pi$.
\item \ If $F_h$, $F_i$, $F_j$, $F_k$ form a 4-prismatic element, then $\alpha_{hi}+\alpha_{ij}+\alpha_{jk}+\alpha_{ki}<2\pi$.
\item \ The assignment of dihedral angles shown in Fig.~\ref{f5.3} does not occur.
\end{enumerate}
Moreover, these conditions are sufficient for an abstract polyhedron $P$ with $\ge 5$ vertices to be realizable as a compact convex polyhedron in
$\h^3$ with dihedral angles equal to $\alpha_{ij}$.}

We shall call such polyhedra \emph{Andreev polyhedra} and conditions 1-5 \emph{Andreev conditions}. Because the realizability is not satisfied when
$P$ is a simplex (an abstract polygon with $4$ vertices), Hodgson restricts himself to the case in which $P$ is not a simplex. The same is assumed
here.

We first see that indeed such a polyhedron has trivalent vertices.

\ite{Lemma}{}{A convex spherical polygon with side lengths $\,\in [\pi/2,\pi)$ is a spherical triangle. If there is a side of length $\pi$ then it is
a bigon.}

\begin{figure}
\begin{center}
\psfrag{l1}{$\ell_1$}\psfrag{l2}{$\ell_2$}\psfrag{lN}{$\ell_N$} \psfrag{r1}{$r_1$}\psfrag{r2}{$r_2$}\psfrag{r3}{$r_3$}\psfrag{rN}{$r_N$}
\psfrag{a1}{$\alpha_1$}\psfrag{a2}{$\alpha_2$}\psfrag{aN}{$\alpha_N$}
\includegraphics[height=7cm]{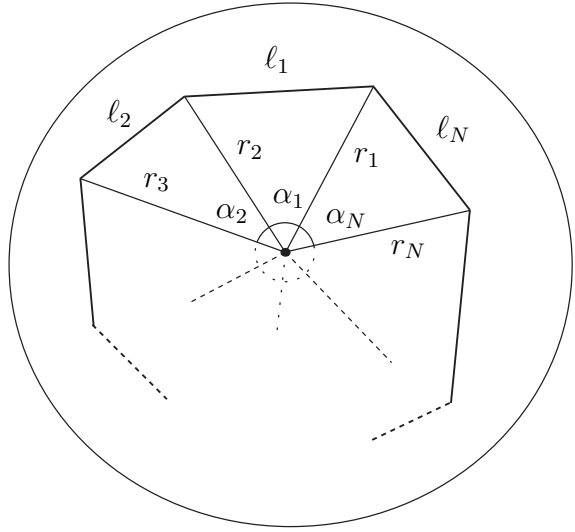}
\end{center}
\caption{Stellation of a spherical polygon.} \label{f5.4}
\end{figure}

\proof Being a convex polygon, it is contained in a hemisphere, which can be chosen with its center $c$ inside the polygon. Let $N$ be the number of
sides. Then it can be triangulated or \emph{stellated} into $N$ spherical triangles, each with a vertex in $c$ and opposite side a side of the
polygon (see Fig.~\ref{f5.4}). The lengths of the sides of the polygon are denoted by $\ell_1,\dots,\ell_N$, the angles at $c$ of the corresponding
triangles are $\alpha_1,\dots,\alpha_N$, and the length of the radii are denoted by $r_1,\dots,r_N$, so that for example the triangle with angle at
$c$ equal to $\alpha_i$ has opposite side of length $\ell_i$ and adjacent sides of length $r_i$ and $r_{i+1}$. By the spherical law of cosines we
have that
\eq{\cos \ell_i = \cos r_i \cos r_{i+1} + \sin r_i \sin r_{i+1} \cos\alpha_i }
The fact that the polygon is contained in the hemisphere implies that $0<r_i<\pi/2$ for all $i$, so $\cos r_i,\sin r_i>0$. Also, by hypothesis
$\ell_i\ge \pi/2$ for all $i$, so $\cos\ell_i\le 0$. Therefore from the spherical law of cosines follows that $\cos\alpha_i<0$ and consequently
$\alpha_i>\pi/2$ for $i=1,\dots,N$. Since $c$ is a smooth point then $\alpha_1+\cdots+\alpha_N=2\pi$, so necessarily $N=3$ unless one of the sides is
of length $\pi$, in which case $N=2$. \qed

The vertices of an Andreev polyhedron correspond under the polar map to convex spherical polygons with sides of length $\pi/2\le\ell_i<\pi$ ($\pi$
minus the corresponding dihedral angle), and by Lemma~5.3 these are necessarily spherical triangles, so Andreev polyhedra have trivalent vertices.

\begin{figure}
\begin{center}
\psfrag{p}{$\pi/2$}
\includegraphics[height=5cm]{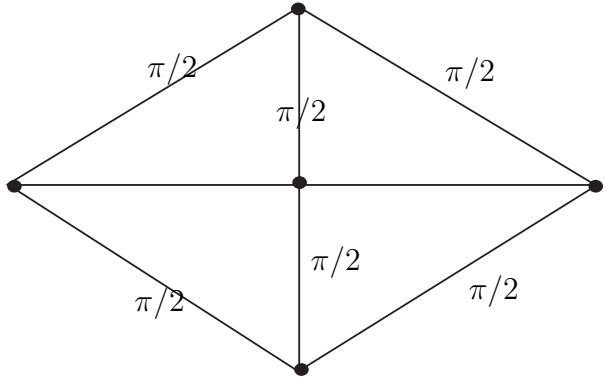}\hspace{1cm}
\end{center}
\caption{A configuration that cannot be part of the polar of a compact hyperbolic polyhedron with all dihedral angles less or equal to $\pi/2$.}
\label{f5.5}
\end{figure}

Now we translate Andreev conditions to conditions on $P^*$ (according to step 1 above) to characterize the set of \emph{Andreev metrics} on $\s^2$.
The \emph{dual Andreev conditions} can be stated as:
\begin{itemize}
\item[$1^*$] $P^*$ can be triangulated by convex spherical triangles with all edge lengths $\ \in[\pi/2,\pi)$.
\item[$2^*$] The lengths of geodesics of $P^*$ comprised of three or four edges of the triangulation is greater than $2\pi$.
\item[$3^*$] $P^*$ does not contain the configuration of Fig.~\ref{f5.5}.
\end{itemize}
Moreover, if $P$ is not a simplex then neither is $P^*$.

Condition $1^*$ follows from two facts. First, edge lengths are complementary to dihedral angles under the polar map, which are in the interval
$(0,\pi/2]$ according to condition $1$. Second, the polar of condition $2$ implies that the perimeter of the triangles of the triangulation of $P^*$
is smaller than $2\pi$, so these triangles are convex.

Condition $3^*$ is just the dual of condition 5. Condition $2^*$ follows from the dual of conditions $3$ and $4$ and from the following lemma:

\ite{Lemma}{(Geodesic edge cycles)}{A closed curve in the 1-skeleton of $P^*$ is a geodesic unless it contains a vertex with a single triangle on one
side.}

\proof First we note a property of convex spherical triangles with edge lengths in the interval $[\pi/2,\pi)$. Let one of such triangles have an
angle $\alpha$ at some vertex, the opposite side of length $a$ and the other two sides $b$ and $c$. By the spherical law of cosines,
$\cos\alpha=(\cos a-\cos b\cos c)/\sin b\sin c$. Since $\sin b\sin c\in (0,1]$ and $(\cos a-\cos b\cos c)\le 0$, we have necessarily that
$\cos\alpha\le(\cos a-\cos b\cos c)$. Also $\cos b\cos c\in(-1,0]$, so $\cos\alpha\le\cos a$, and since the cosine is decreasing in $[\pi/2,\pi]$, it
follows that $\alpha\ge a$. So each angle is no less than the opposite side. This proves that a curve through a cone point in $P^*$ that leaves at
least two triangles on each side subtends an angle $\ge\pi$ at each side. \qed

Now we follow step 2 and prove that the polar of an Adreev polyhedron $P$ is an Andreev metric, that is, that satisfies the dual Andreev conditions.
Let $P^*$ be the dual of $P$. Since $P$ has trivalent vertices, $P^*$ is the topological sphere obtained by gluing spherical triangles, with side
lengths which are dual of dihedral angles of $P$, and therefore in the interval $[\pi/2,\pi)$. This gives the desired triangulation (moreover this
triangulation is unique). Since $P^*$ defines an admissible metric according to Rivin's theorem, then the length of geodesics of $P^*$ is always
greater than $2\pi$, so $P^*$ satisfies $2^*$. In addition, $P^*$ cannot contain the configuration of Fig.~\ref{f5.5}, since the cone angle at the
interior vertex is necessarily $2\pi$, contradicting Rivin's theorem.

Step 3 requires to prove that Andreev metrics are admissible. This follows from two propositions:

\ite{Proposition}{}{A metric space $P^*$ satisfying the dual Andreev conditions has all cone angles greater than $2\pi$.}

\proof Consider the cone angle at the vertex $v$, and the star of $v$, denoted by $st(v)$. Since, as shown before, an angle of a spherical triangle is greater or equal to the length of the opposite side, we must only check the cases in which the star of $v$ is composed of 2 or 3 triangles, or of 4 triangles with opposite sides of length $\pi/2$. This last case is ruled out by condition $3^*$. The case in which only two triangles are incident on $v$ does not
occur in a valid triangulation. Suppose there are exactly 3 triangles incident on $v$. Then the closed curve composed of the three opposite sides of
those triangles is a geodesic (if it passed through a vertex with one triangle on the other side then $P^*$ would be a simplex). So by condition
$2^*$ its length is greater than $2\pi$. The cone angle at $v$ is greater or equal to the length of the curve around the star, so it is greater than
$2\pi$. \qed

\ite{Proposition}{}{If $P^*$ is a metric space satisfying the dual Andreev conditions, then every closed geodesic in $P^*$ has length $>2\pi$.}

\proof For this proof we need two preliminary facts. Let $v$ be a vertex of $P^*$ with cone angle $>2\pi$, and let $ost(v)$ and $\partial st(v)$ be the interior and the border of $st(v)$. Then,

\begin{enumerate}
\item Every geodesic contained in  $ost(v)$ is strictly longer than $2\pi$.
\item A geodesic arc $\gamma$ intersecting $ost(v)$ and joining two points $p$ and $q$ in $\partial st(v)$ has length greater or equal to $\pi$. If its length is exactly $\pi$, then both $p$ and $q$ are at distance $\pi/2$ from $v$ and either $\gamma$ passes through $v$ or there is another geodesic of length $\pi$ in $\partial st(v)$ such that both geodesics bound a lune not containing $v$. 
\end{enumerate}

Now consider a geodesic $\gamma$ in $P^*$. If $\gamma$ is contained in the star of a vertex, then it has length $>2\pi$. If $\gamma$ intersects three or more disjoint stars, then it has length $\ge 3\pi$. Therefore we must only consider the case in which $\gamma$ is contained in the union of the stars of two disconnected vertices $v_1$ and $v_2$. We call $\gamma_1=\gamma\cap st(v_1)$ and $\gamma_2=\gamma\cap st(v_2)$. We must prove that it is not possible that $\gamma$ is of length $2\pi$, which by fact~2 above is equivalent to proving that it is not possible that both $\gamma_1$ and $\gamma_2$
have length $\pi$.

 Suppose that both $\gamma_1$ and $\gamma_2$ have length $\pi$. Then their endpoints $p$ and $q$ are at distance $\pi/2$ from $v_1$ and $v_2$.  There are several possibilities:
 
 \begin{enumerate}
 \item If both $v_1$ and $v_2$ lie on $\gamma$, then the segments $v_1p$, $v_1q$ etc., can be replaced by segments of the triangulation such that $\gamma$ is replaced by a geodesic $\gamma'$ of length $2\pi$ in the 1-skeleton of $P^*$. The fact that $\gamma'$ is a geodesic is guaranteed by condition $3^*$.
 \item  If both $p$ and $q$ are vertices of $P^*$ then the previous argument applies replacing $v_1$ and $v_2$ with $p$ and $q$, so again there is a geodesic $\gamma'$ of length $2\pi$ in the 1-skeleton of $P^*$.
 \item If $\gamma_1$ does not contain $v_1$ and $p$ is not a vertex of $P^*$, then there is a geodesic arc $\gamma_1'$  in the 1-skeleton of $P^*$ of length $\pi$ from $p$ to $q$. Moreover $\gamma_1$ and $\gamma_1'$ bound a lune of angle $<\pi/2$.  Since $\gamma_1$ and $\gamma_2$ join at angle $\pi$ at $p$, it follows that $\gamma_2$ cannot contain $v_2$, and there is another geodesic segment $\gamma_2'$ of length $\pi$ in the 1-skeleton of $P^*$ joining $p$ and $q$ such that $\gamma_2$ and $\gamma_2'$ bound a lune of angle $<\pi/2$. By angle considerations it can be seen that $\gamma_1'\ne\gamma_2'$, so $\gamma'=\gamma_1'\cup\gamma_2'$ is a geodesic of length $2\pi$ in the 1-skeleton of $P^*$.
 \end{enumerate}
 
 But by condition $2^*$ a geodesic of length $2\pi$ cannot be contained in the 1-skeleton of $P^*$, so the possibility that both $\gamma_1$ and $\gamma_2$ have length $\pi$ is ruled out. \qed

\section{Pogolerov's map and some counter-examples}

\setcounter{sec}{1}

The purpose of this section is to prove that edge lengths do not determine hyperbolic polyhedra even when the combinatorics is fixed. The proof consists on the construction of explicit counter-examples, as given by Schlenker in \cite{Schlenker}. Since the main tool used to construct such counter-examples is the \emph{Pogolerov's map}, we shall give a detailed description of this application before. 

We begin writing the explicit form for the (projective) application that takes points in the hyperboloid model of $\h^3$ to the
projective model. This is the map $\rho:\h^3\to B^3$ given by
\eq{\rho(x)=\frac{\vec{x}}{x_0}}
where $x=(x_0,\vec{x})$ is any point on the hyperboloid.

This map can be generalized in an ingenious way through an application due to Pogorelov \cite{Pogol,Schlenker}. We will first define Pogorelov's map
and present some properties. It will then become clear in what sense it is a generalization of the projective application $\rho$. Contrary to what is
done in the cited references, here we choose to present constructive proofs, giving explicit expressions for inverse mappings and isometries in terms
of coordinates.

\ite{Definition}{}{Let $x=(x_0,\vec{x})$ and $y=(y_0,\vec{y})$ be two arbitrary points in the hyperboloid model of $\h^3$. The {\bf Pogorelov's map},
$\Phi$, is defined as the mapping $\Phi:\h^3\times\h^3\to\R^3\times\R^3$ such that
\eq{\Phi(x,y)=\Bigg( \frac{2\vec{x}}{x_0+y_0},\frac{2\vec{y}}{x_0+y_0} \Bigg)}}

This map is continuous, since $x_0+y_0\ge 2$ for points on the hyperboloid. It is also injective, as can be checked by direct computation:

\ite{Lemma}{}{Pogolerovs's map is a homeomorphism onto its image.}

\proof By the theorem of invariance of domain, since $\Phi$ is continuous, it remains to prove that $\Phi$ is one-to-one. Let $x,y,x',y'\in \h^3$ such that $\Phi(x,y)=\Phi(x',y')$. Then
\eq{\frac{2\vec{x}}{x_0+y_0}=\frac{2\vec{x}^{\,\prime}}{x'_0+y'_0}\quad,\qquad \frac{2\vec{y}}{x_0+y_0}=\frac{2\vec{y}^{\,\prime}}{x'_0+y'_0}\ ,}
so we can write $\vec{x}^{\,\prime}=\lambda\vec{x}$ and $\vec{y}^{\,\prime}=\lambda\vec{y}$ with $\lambda=(x_0'+y_0')/(x_0+y_0)>0$. Now, since
$x,x',y$ and $y'$ are on the hyperboloid, we have
\eqa{x^{\prime\,2}_0=1+\|\vec{x}^{\,\prime}\|^2=1+\lambda^2(x_0^2-1)\nonumber\\[2mm]
y^{\prime\,2}_0=1+\|\vec{y}^{\,\prime}\|^2=1+\lambda^2(y_0^2-1)\nonumber}
and substituting this into the expression for $\lambda$ leads to the equation
\eq{\lambda\,(x_0+y_0)=\sqrt{1+\lambda^2(x_0^2-1)}+\sqrt{1+\lambda^2(y_0^2-1)}\ .}
Solving this equation for $\lambda$ gives $\lambda=\pm 1$, and since $\lambda>0$ then $\lambda=1$. So $\vec{x}^{\,\prime}=\vec{x}$ and
$\vec{y}^{\,\prime}=\vec{y}$, which implies that $x=x'$ and $y=y'$, because all four points lay on the hyperboloid. \qed

Therefore Pogolerov's map is a homeomorphism onto its image (in fact it is diffeomorphism, but we don't need that). The image is
\eq{{\rm Im}(\Phi)=\{(a,b)\in \R^3\times\R^3\,|\,\|a\|+\|b\|<2\}\ ,}
although we will not prove this. We shall prove however a weaker result that will be useful later:

\ite{Proposition}{}{$B^3\times B^3\subset {\rm Im}(\Phi)$.}

\proof Let $\vec{a},\vec{b}\in B^3$ be two arbitrary vectors in the unit 3-ball. They can be expressed as $\vec{a}\equiv a\hat{u}$ and $\vec{b}\equiv
b\hat{v}$, where $\hat{u},\hat{v}$ are unit vectors and $a,b\in[0,1)$. Choose the following two points $x=(x_0,\vec{x})$, $y=(y_0,\vec{y})$ on the
hyperboloid:
\eqa{ \vec{x} &=& \frac{4a\,\hat{u}}{\sqrt{f(a,b)}}\quad ,\quad x_0=\sqrt{1+\|\vec{x}\|^2}\nonumber\\[2mm]
 \vec{y} &=& \frac{4b\,\hat{v}}{\sqrt{f(a,b)}}\quad ,\quad y_0=\sqrt{1+\|\vec{y}\|^2}\nonumber  }
with $f(a,b)=(a^2-b^2)^2-8(a^2+b^2-2)$. These points are well defined for all $a,b\in[0,1)$, because $f(0,0)=16$, $f$ is continuous, and $f\ne 0$ for
all $a,b\in[0,1)$, so $f>0$ in that region. Now it is straightforward to verify that $\Phi(x,y)=(\vec{a},\vec{b})$, which proves that
$(\vec{a},\vec{b})\in {\rm Im}(\Phi)$. \qed

The \emph{diagonal} subset of $\h^3\times\h^3$ is the set $\Delta=\{(x,x)\,|\,x\in\h^3\}\subset \h^3\times\h^3$. When restricted to $\Delta$, the
action of $\Phi$ is familiar: $\Phi(x,x)=(2\vec{x}/2x_0,2\vec{x}/2x_0)=(\rho(x),\rho(x))$. This simple fact is the content of the following
proposition.

\ite{Proposition}{}{The restriction of $\Phi$ to the diagonal corresponds to the projective application $\rho$ on each factor.}

Now we prove the two more interesting properties of the Pogolerov's map. First, consider two congruent figures in $\h^3$, $F$ and $F'$. Since they
are congruent, one can be obtained from the other by applying an isometry $A$, $F'=AF$. (Note that as we are working in the hyperboloid model of
$\h^3$, $A$ corresponds to a (proper) Lorentz transformation.) We apply the Pogorelov's map to the pair $(F,AF)$ in the sense that
$\Phi(F,AF)=\{\Phi(x,Ax)\,|\,x\in F\}$. This gives $\Phi(F,AF)=(G,G')$, where $G,G'$ are figures in $\R^3$. The first nice property of $\Phi$ is that
$G$ and $G'$ are congruent, and the Euclidean isometry that relates them depends exclusively on the hyperbolic isometry $A$.

\ite{Proposition}{}{Let $A$ denote a hyperbolic isometry, and $x\in \h^3$. We write $\Phi(x,Ax)=(y,y')$ with $y,y'\in \R^3$. Then, for each $A$ there
is an Euclidean isometry $B$ such that, for all $x$, $y'=By$.}

\proof Let $\vec{a},\vec{b}\in \R^3$ such that $(\vec{a},\vec{b})\in {\rm Im}(\Phi)$, and let $a=\|\vec{a}\|$ and $b=\|\vec{b}\|$. The inverse of
$(\vec{a},\vec{b})$ under $\Phi$ is given by $\Phi^{-1}(\vec{a},\vec{b})=(x,y)$ with (see the proof of Proposition 6.3)
\eqa{x=\frac{1}{\sqrt{f(a,b)}}\,\Big(4+a^2-b^2\,,\,4\vec{a}\,\Big)\nonumber\\[2mm]
y=\frac{1}{\sqrt{f(a,b)}}\,\Big(4-a^2+b^2\,,\,4\vec{b}\,\Big)\nonumber}
Now we consider the case in which $y=Ax$, with $A\in {\rm Isom}(\h^3)$. This implies the following equations for $\vec{a}=(a_1,a_2,a_3)$ and
$\vec{b}=(b_1,b_2,b_3)$:
\eqa{(4-a^2+b^2)&=&A_{00}\,(4+a^2-b^2)+4A_{0j}\,a_j\nonumber\\[2mm]
4\,b_i &=& A_{i0}\,(4+a^2-b^2)+4A_{ij}\,a_j\nonumber}
where a summation for $j=1,2,3$ is understood whenever the index $j$ appears twice, and $A_{kl}$ are the components of $A$ in the basis of $\e^3_1$.
From these equations the following relationship between $\vec{a}$ and $\vec{b}$ arises,
\eq{\vec{b}=\vec{D}+R\vec{a}\ ,}
where $\vec{D}$ is a 3-displacement and $R$ a 3-rotation given in components by
\eq{D_i=\frac{2A_{i0}}{1+A_{00}}\quad,\qquad R_{ij}=A_{ij}-\frac{A_{i0}A_{0j}}{1+A_{00}}\ .}
Using the fact that $A$ is a Lorentz transformation (more specifically, using the identity $A_{ij}A_{kj}=\delta_{ik}+A_{i0}A_{k0}$), it can be
checked explicitly that $R_{ij}R_{kj}=\delta_{ik}$, so $R$ is indeed a 3-rotation. Note also that $A_{00}=-1$ does not occur in an isometry of
$\e^3_1$ that fixes the positive (or negative) hyperboloid, where necessarily $A_{00}\ge 1$. Then the Euclidean isometry $B$ such that
$B\vec{v}=\vec{D}+R\vec{v}\ $ for all $\vec{v}\in\R^3$, is the one we were looking for. \qed

The second interesting property of $\Phi$ is that pairs of geodesic planes in $\h^3$ are mapped to pairs of geodesic planes in $\R^3$.

\ite{Proposition}{}{$\Phi$ is a geodesic mapping.}

\proof Let $p_1$ be the projection of $\R^3\times\R^3$ on the first factor, and $\vec{\Phi}_1=p_1\circ\Phi$. Let $A\in{\rm Isom}(\h^3)$. Consider all
$x\in\h^3$ satisfying the following linear equation,
\eq{\vec{a}\cdot\vec{\Phi}_1(x,Ax)+b=0}
with $\vec{a}$ a unit three vector and $b$ a scalar. We denote this set by $X$, and it is the set in $\h^3$ such that $P=p_1\circ\Phi(X,AX)$ defines
the 2-plane in $\R^3$ perpendicular to $\vec{a}$ and a distance $b$ from the origin. Using the explicit form of $\vec{\Phi}_1$ in coordinates, $X$ is
the set of points $x=(x_0,\vec{x})\in\e^3_1$ such that
\eqa{1.&&2\,\vec{a}\cdot \vec{x}+b\,[x_0+(Ax)_0]=0\nonumber\\[2mm]
2.&&x\cdot x=-1\nonumber}
The fist condition is linear and homogeneous, so it defines a plane in $\e^3_1$ through the origin. The second condition means that we consider the
points of intersection of such a plane with the hyperboloid. So $X$ is precisely a 3-plane in $\h^3$, as well as $AX$. Moreover, from Proposition 6.5
we know that there is an Euclidean isometry $B$ such that $\Phi(X,AX)=(P,BP)=(P,P')$, with $P'$ a 2-plane in $\R^3$.

Therefore, $\Phi$ takes pairs of planes to pairs of planes. This also means that it takes pairs of geodesics to pairs of geodesics, since each
geodesic is the intersection of two non-parallel planes. \qed

To finish with the set of properties of Pogolerov's map, there is one last result that we will need.

\ite{Proposition}{}{Let $P$ be a plane in $\R^3$ containing the origin. Then the set of points $x\in\h^3$ such that $\vec{\Phi}_1(x,y)\in P$ for all
$y\in\h^3$, is the plane $\rho^{-1}(P)$. Conversely, for any plane $H$ in $\h^3$ containing $\rho^{-1}(0)$, the set of points $\vec{a}\in\R^3$ such
that $\Phi^{-1}(\vec{a},\vec{b})\in H\times\h^3$ for all $\vec{b}\in {\rm Im}(\vec{\Phi}_2)$, is a plane in ${\rm Im}(\Phi_1)$ containing $\rho(H)$.}

\proof A plane $P$ though the origin in $\R^3$ is the set $P=\{\vec{a}\in\R^3\,|\,\vec{v}\cdot\vec{a}=0\}$ for a fixed $\vec{v}\in\R^3$. Now consider
the set of points $X=\{x\in\h^3\,|\,\vec{\Phi}_1(x,y)\in P\ \forall\,y\in\h^3\}$. Then, since $\vec{\Phi}_1(x,y)=\lambda \vec{x}$ with $\lambda\ne
0$, we have that $\vec{\Phi}_1(x,y)\in P\Leftrightarrow \vec{x}\in P$, so
\eqa{X&=&\{x\in\h^3\,|\,\vec{x}\in P\}=\{x\in\h^3\,|\,\vec{v}\cdot\vec{x}=0\}=\{x\in\h^3\,|\,\vec{v}\cdot (\vec{x}/x_0)=0\}\nonumber\\[2mm]
&=&\{x\in\h^3\,|\,\vec{v}\cdot \rho(x)=0\}=\{x\in\h^3\,|\,\rho(x)\in P\}=\rho^{-1}(P)\ .\nonumber}

Now, a plane $H$ in $\h^3$ containing $\rho^{-1}(0)$ is the intersection with the hyperboloid of a plane
$\widetilde{H}=\{z\in\e^3_1\,|\,\vec{z}\cdot\vec{v}=0\}\subset\e^3_1$, for a fixed $\vec{v}\in\R^3$. Consider the set of points $Y=\{\vec{a}\in{\rm
Im}(\vec{\Phi}_1)\,|\,\Phi^{-1}(\vec{a},\vec{b})\in H\times\h^3\ \forall\,\vec{b}\in{\rm Im}(\vec{\Phi}_2)\}$. From the explicit expression for
$\Phi^{-1}$ we see that $\Phi^{-1}(\vec{a},\vec{b})\in H\times\h^3\Leftrightarrow \vec{a}\cdot\vec{v}=0$. Also, for $\vec{a}\in B^3$ we have that
$\rho^{-1}(\vec{a})\propto (1,\vec{a})$, so
\eqa{Y&=&\{\vec{a}\in{\rm Im}(\vec{\Phi}_1)\,|\,\vec{a}\cdot\vec{v}=0\}=\widetilde{H}\cap {\rm Im}(\Phi_1)\supset\{\vec{a}\in
B^3\,|\,\vec{a}\cdot\vec{v}=0\}\nonumber\\[2mm]
&=&\{\vec{a}\in B^3\,|\,(1,\vec{a})\in \widetilde{H}\}=\{\vec{a}\in B^3\,|\rho^{-1}(\vec{a})\in \widetilde{H}\}=\rho(H)\ .\nonumber}
Clearly $Y\ne\rho(H)$ because $Y$ contains points outside the unit 3-ball. \qed

Having introduced Pogolerov's map and its properties, we can now use it to construct examples of pairs of non-congruent polyhedra in $\h^3$ with the
\emph{same} combinatorics and the \emph{same} edge lengths \cite{Schlenker}. Note that from Rivin's theorem does not follow that edge lengths
determine polyhedra in $\h^3$, but it is also not easy to prove the contrary. The counter-example in Ref.~\cite{Schlenker} proves this.

The idea to construct such counter-examples is the following: edge lengths do not determine Euclidean polyhedra, so we can construct pairs of
non-congruent Euclidean polyhedra inside $B^3$, and then use the inverse of Pogolerov's map to obtain pairs of polyhedra in $\h^3$. The properties of
this map ensure (if the Euclidean pairs are chosen correctly) that these pairs of hyperbolic polyhedra are non-congruent and have the same edge
lengths.

\begin{figure}
\begin{center}
\psfrag{a}{$a$}\psfrag{b}{$b$}\psfrag{c}{$c$}\psfrag{u}{$u$}\psfrag{x}{$x$}\psfrag{y}{$y$}\psfrag{z}{$z$}
\includegraphics[height=7cm]{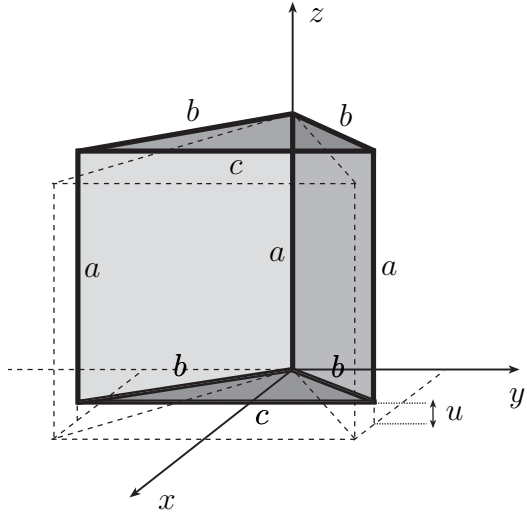}
\end{center}
\vspace{-0.5cm} \caption{The prism $P_{(a,b,c)}^{\,u}$ compared with the prism $P_{(a,b,c)}^{\,0}$ (dashed).} \label{prisma}
\end{figure}

Consider an abstract triangular prism $P$. We first construct a family of polyhedral immersions $\alpha_u:P\to B^3$, such that the polyhedra
$\alpha_u(P)$ are all Euclidean triangular prisms with the same edge lengths, but non-congruent for different values of $u$. These immersions are
defined such that $\alpha_u(P)=P_{(a,b,c)}^{\,u}$, where $P_{(a,b,c)}^{\,u}$ is the prism with vertices $(0,0,0)$, $(0,0,a)$,
$(\sqrt{b^2-c^2/4-u^2},c/2,u)$, $(\sqrt{b^2-c^2/4-u^2},-c/2,u)$, $(\sqrt{b^2-c^2/4-u^2},c/2,u+a)$ and $(\sqrt{b^2-c^2/4-u^2},-c/2,u+a)$, as shown
in Fig.~\ref{prisma}. Clearly, for $u<\sqrt{b^2-c^2/4}$ these prisms are convex and have the same edge lengths: $a$, $b$ and $c$. Moreover, for
$a,b,c$ small enough $P_{(a,b,c)}^{\,u}\subset B^3$. The family of immersions $\{\alpha_u\}$ is then defined for $u<\sqrt{b^2-c^2/4}$ and for small
and fixed $(a,b,c)$.

Now we consider two prisms $\alpha_u(P)$ and $\alpha_v(P)$ with $u\ne v$. By definition $(\alpha_u,\alpha_v)(P)\in B^3\times B^3$, so by
Proposition~6.3, $(\alpha_u,\alpha_v)(P)\in{\rm Im}(\Phi)$. Then by Proposition~6.2 we know that $\Phi^{-1}$ is a well defined application, and we
can consider the image in $\h^3\times\h^3$ of $P$ under $\Phi^{-1}\circ (\alpha_u,\alpha_v)$. We want to prove that for $(F,F')=\Phi^{-1}\circ
(\alpha_u,\alpha_v)(P)$, $F$ and $F'$ are convex non-congruent prisms in $\h^3$ with the same edge lengths.

That $F$ and $F'$ have the same edge lengths follows from Proposition~6.5 and Proposition~6.6, since pairs of congruent edges of $\alpha_u(P)$ and
$\alpha_v(P)$ are mapped under $\Phi^{-1}$ to pairs or congruent edges of $F$ and $F'$. Also, for the same reason, the two pairs of triangular faces
and the pair of rectangular faces of side lengths $a$ and $c$ (with are congruent in $\alpha_u(P)$ and $\alpha_v(P)$) are mapped to congruent faces
in $F$ and $F'$. Note that it is not only important that these edges and faces are congruent in $F$ and $F'$, but also that correspond to geodesic
edges and faces, since it must be proven that $F$ and $F'$ are polyhedra.

The faces of side lengths $a$ and $b$ are not congruent in $\alpha_u(P)$ and $\alpha_v(P)$, so the previous argument does not apply. However they
both contain the origin, so Proposition~6.7 tells us that they are mapped inside hyperbolic planes. The fact that they do indeed constitute the faces
that are missing so far in $F$ and $F'$ follows from continuity and taking a pair of congruent vertical segments inside these faces. So $F$ and $F$
are hyperbolic polyhedra with the same edge lengths. Since $\rho$ is projective and $\Phi$ restricted to the diagonal is $\rho$ on each factor (by
Proposition~6.4), it follows from continuity that for $(u-v)$ small enough, $F$ and $F'$ are convex.

The fact that $F$ and $F'$ are non-congruent follows from Proposition~6.5, because if they were congruent then $\alpha_u(P)$ and $\alpha_v(P)$ would
be congruent, which is not true by hypothesis. So $F$ and $F'$ are convex non-congruent prisms in $\h^3$ with the same edge lengths. This proves the
following theorem:

\ite{Theorem}{}{There exist pairs of non-congruent polyhedra in $\h^3$ with the same edge lengths.}

It would probably be interesting to understand these counter-examples in the context of Rivin's characterization, that is, to see how these families of polyhedra look like in the space of admissible metrics. It could also help to establish a systematic procedure to build families of non-congruent convex hyperbolic polyhedra with the same combinatorics and edge lengths, or even to find the complete set of such families, but maybe this is too optimistic.



%
%
%
%
%
%
%
%
%
%


\begin{thebibliography}{99}

\bibitem{Rivin1}
   I.~Rivin and C.~D.~Hodgson,
   \emph{A characterization of compact convex polyhedra in hyperbolic 3-space}. Invent.Math., 111:77-111, 1993.

\bibitem{Rivin4}
   I.~Rivin,
   \emph{On geometry of convex ideal polyhedra in hyperbolic 3-space}. Topology 32:87-92, 1993.

\bibitem{Rivin2}
   I.~Rivin,
   \emph{A characterization of ideal polyhedra in hyperbolic 3-space}. Annals of Mathematics, 143:51-70, 1996.

\bibitem{Hodgson}
   C.~D.~Hodgson,
   \emph{Deduction of Andeev's theorem from Rivin's characterization of convex hyperbolic polyhedra}.
   Topology 90. Proceedings of the Research Semester in Low Dimensional Topology at O.S.U Berlin New York: de Gruyter, 1993.

\bibitem{Schlenker}
   J.~M.~Schlenker,
   \emph{Dihedral angles of convex polyhedra}. Discrete Comput. Geom. 23:409-417, 2000.

\bibitem{Rivin3}
   C.~D.~Hodgson, I.~Rivin and W.~D.~Smith,
   \emph{A characterization of convex hyperbolic polyhedra and of convex polyhedra inscribed in the sphere}.
   Bulletin of the American Mathematical Society, 27:246-251, 1992.

\bibitem{Bobenko}
   A.~Bobenko and I.~Izmestiev,
   \emph{Alexandrov's theorem, weighted Delaunay triangulations, and mixed volumes}. arXiv:math/0503219, 2006.

\bibitem{Andreev}
   E.~M.~Andreev,
   \emph{On convex polyhedra in Lobachevskii space}. Math USSR. Sbornik. 10:413-440, 1970.

\bibitem{Roeder1}
   R.~K.~W.~Roeder, J.~H.~Hubbard and W.~D.~Dunbar,
   \emph{Andreev's Theorem on  hyperbolic polyhedra}. Annales de l'institut Fourier, 57 no.~3, p.~825-882, 2007.

\bibitem{Pogol}
   A.~V.~Pogorelov,
   \emph{Extrinsic Geometry of Convex Surfaces}. Translations of Mathematical Monographs, Vol.35, American Mathematical Society, 1973.

\bibitem{Th}
   W.~Thurston,
   \emph{Three-dimensional geometry and topology}. Princeton University Press, Princeton, NJ, 1997. Edited by Silvio Levy.

\bibitem{Rat}
   J.~G.~Ratcliffe,
   \emph{Foundations of hyperbolic manifolds}. Graduate Texts in Mathematics, 149. Springer-Verlag, New York, 1994.




\end{thebibliography}
\end{document}